\documentclass{article}
\usepackage[utf8]{inputenc}
\usepackage{amssymb, amsmath, amsfonts, blkarray, amsthm}
\usepackage[round]{natbib}
\usepackage{hyperref}
\usepackage{enumerate}
\usepackage{titlesec}
\usepackage{parskip}
\usepackage{listings} 
\usepackage[normalem]{ulem} 
\usepackage{xcolor, graphicx}
\usepackage{tabularx}
\usepackage{stmaryrd}
\usepackage{standalone}
\usepackage{tikz}
\usepackage{authblk}
\usetikzlibrary{arrows.meta,calc,decorations.markings,angles,quotes}
\usepackage{pgfplots}
\pgfplotsset{compat=1.18}

\newtheorem{theorem}{Theorem}[section]
\newtheorem{corollary}[theorem]{Corollary}
\newtheorem{lemma}[theorem]{Lemma}
\newtheorem{prop}[theorem]{Proposition}
\newtheorem{remark}[theorem]{Remark}

\begin{document}
\title{Group-averaged Markov chains II: tuning of group action in finite state space}

\author[1]{Michael C.H. Choi\thanks{Email: mchchoi@nus.edu.sg, corresponding author}}

\author[1]{Ryan J.Y. Lim\thanks{Email: ryan.limjy@u.nus.edu}}

\author[1]{Youjia Wang\thanks{Email: e1124868@u.nus.edu}}

\affil[1]{Department of Statistics and Data Science, National University of Singapore, Level 7, 6 Science Drive 2, 117546, Singapore}
\date{\today}

\maketitle

\begin{abstract}
We study group-averaged Markov chains obtained by augmenting a $\pi$-stationary kernel $P$ with orbit kernels induced by a group action. We analyse the Gibbs ($G$), Metropolis--Hastings ($M$), and Barker ($B$) kernels, their sandwiches $QPQ$, and mixtures $\tfrac{1}{2}(P+Q)$, where $Q\in\{G,M,B\}$. Under suitable conditions, $M^t$ and $B^t$ converge blockwise to $G$. The projection chains of $GPG$ and $P$ coincide, while every sandwich $QPQ$ has absolute spectral gap no smaller than that of reversible $P$. For $GPG$, we derive an additive asymptotic-variance bound, prove monotonicity for $G$-invariant observables, and identify it as the Kullback--Leibler (KL) information projection of $P$ onto the $G$-invariant kernels. For a fixed orbit partition, the spectral and KL properties of $GPG$ reduce to those of a lower-dimensional orbit-space chain. Among Gibbs projections with a prescribed number of orbits, we identify the partition minimizing KL divergence to stationarity and characterize exact stationarity. Finally, alternating group projections converge at a rate determined by singular values of an overlap matrix and, in structured cases, can yield exact sampling with logarithmically many group actions. These results motivate tuning heuristics and yield polynomial mixing for a Curie--Weiss example in a regime where Glauber dynamics is exponentially slow.\\

\end{abstract}

\tableofcontents

\section{Introduction}
While Markov chain Monte Carlo (MCMC) methods remain indispensable for sampling from complex and high-dimensional distributions, their efficiency often deteriorates when the target distribution exhibits strong multimodality. In such settings, standard chains easily become trapped within local regions of the state space. Recent work has explored ways to accelerate mixing by augmenting Markov chains with structured or deterministic transitions, such as the deterministic jump framework of \cite{Diaconis_2020}. \cite{ying2022doubleflipisingmodels} further introduced a double-flip move for Ising models, implemented as an additive mixture on top of the Swendsen-Wang algorithm. This construction can be interpreted as a special case of a symmetry-based jump, equipped with a Metropolis–Hastings correction, generated by the two-element cyclic group $C_2$.

Building on the approaches introduced in \cite{Choi_Permutation} and \cite{choi2025groupaveragedmarkovchainsmixing}, this paper develops a systematic method for incorporating group actions to improve sampling dynamics. Unlike previous formulations that rely on equi-probability jumps, our construction allows general $\pi$ weighted transitions within group orbits, yielding a broader and more flexible class of group-augmented samplers.

Our work fits within a growing line of research that studies Markov chains whose behaviour is shaped by group actions and the orbit partitions they induce. The closest examples are the Burnside processes and its recent developments \cite{Jerrum_1993, aldous-fill-2014, diaconis2025curiouslyslowlymixingmarkov, Diaconis_2025, Diaconis_2021_hahn, diaconis2025countingnumbergrouporbits}. These chains move between group orbits in order to sample uniformly from orbit space, and it demonstrate how orbit structure can produce strikingly different mixing rates. They also support a range of applications, including the simulation of contingency tables and partition-like objects.

The idea of exploiting symmetries and group actions extends well beyond Markov chains. For instance, group equivariant neural networks in \cite{cohen2016groupequivariantconvolutionalnetworks, kondor2018generalizationequivarianceconvolutionneural}, incorporate rotations, reflections, and translations directly into their architecture to enforce invariance and reduce sample complexity. In probabilistic graphical models, the study of automorphism groups by \cite{bui2012automorphismgroupsgraphicalmodels} shows how structural symmetries can be leveraged during inference. Recent advances in generative modelling, including structure-preserving GANs introduced by \cite{birrell2022structurepreservinggans} and group-invariant GANs analysed by \cite{chen2025statisticalguaranteesgroupinvariantgans}, further demonstrate the benefits of embedding symmetry into the model design to enhance accuracy and data efficiency.  These works demonstrate how symmetry can be introduced deliberately to improve efficiency in various contexts.

The main results of this paper can be organised into four parts. We first formalise the construction of several canonical orbit kernels, namely, the Gibbs $(G)$, Metropolis–Hastings $(M)$, and Barker $(B)$ kernels. We then analyse their interactions with a base sampler through multiplicative sandwiches $(QPQ)$ and additive mixtures $\tfrac{1}{2}(P+Q)$ for any $Q \in \{G,M,B\}$. Next, we establish that multiplicative group averaging does no worse than the original sampler in terms of absolute spectral gap and Kullback--Leibler divergence. For asymptotic variance, we derive an additive upper bound for $GPG$ and show, in particular, that the asymptotic variance does not increase for $G$-invariant observables. Further, we show that $GPG$ arises naturally as the information projection of $P$ onto the set of $G$-invariant transition matrices. We then investigate optimality conditions, characterising the optimal sampler $P$ for a given group action $\mathcal{G}$, as well as the optimal $G$ for a fixed $P$. Finally, we explore alternating projections, where multiple group actions ${G_i}$ are composed to form higher-order group-averaged samplers. We further demonstrate that under suitable symmetry or uniformity conditions, such constructions can achieve exact sampling from $\pi$ using only a logarithmic number of group actions relative to the size of the state space.

In further detail, Section \ref{section:prelim} formulates the three canonical orbit kernels, $G$, $M$, and $B$, together with their multiplicative sandwiches $QPQ$ and additive mixtures $\tfrac{1}{2}(P+Q)$. We establish their connections to the projection and restriction chains induced by the group action $\mathcal{G}$, providing a structural interpretation of how group averaging modifies the base dynamics.

Sections \ref{section:specgap} and \ref{section:var} analyse these samplers in terms of absolute spectral gap and asymptotic variance, respectively. We show that each multiplicative sandwich has absolute spectral gap no smaller than that of the original sampler. For asymptotic variance, we show that $GPG$ does not increase the variance of $G$-invariant observables and derive an additive upper bound for general observables. For $GPG$, we further derive a closed-form expression for the absolute spectral gap as a function of the stationary distribution $\pi$.

In Section \ref{section:KL}, we prove that $GPG$ is the exact information projection of $P$ onto the set of $G$-invariant transition matrices, while the Metropolis–Hastings and Barker orbit kernels act as KL-contractive updates converging towards this invariant set. Furthermore, we show that under most conditions, the invariant sets corresponding to the multiplicative sandwiches of $G$, $M$, and $B$ coincide.

Section \ref{section:bestP} shows that, for a fixed orbit partition, the spectrum and KL divergence of $GPG$ reduce exactly to the corresponding quantities of the projection chain on the orbit space. This motivates a lower-dimensional construction of candidate samplers, including a star-shaped orbit-space kernel that becomes increasingly effective when the stationary mass concentrates on one orbit.

Section \ref{section:bestG} instead optimises the Gibbs projection $G$ itself: among partitions with a fixed number of orbits, we identify the partition minimising $D_{KL}^\pi(G\|\Pi)$, where $\Pi$ is the matrix with all rows equal to $\pi$. For this choice, we further characterise the base kernels $P$ for which $GPG=\Pi$.

Next, Section \ref{section:altproj} introduces the framework of alternating projections involving multiple group actions ${G_i}$. We show that the rate of convergence can be characterised by the singular values of a matrix encoding the overlaps between orbit partitions, and that in certain structured cases, this construction yields an exact sampler. The section also demonstrates how the limiting kernel can be determined deterministically from the combined structure of the group actions.

Lastly, Section \ref{section:tuning} concludes the paper by suggesting several heuristics for tuning and selecting appropriate group actions, particularly in settings where no obvious symmetry or relational structure exists in the state space.

\section{Preliminaries} \label{section:prelim}
Let $\mathcal{X}$ be a finite state space, and $\mathcal{P}(\mathcal{X})$ be the set of all probability masses with full support on $\mathcal{X}$. That is, $\min_x\pi(x) > 0$ for all $\pi \in \mathcal{P}(\mathcal{X})$. For integers $a \leq b \in \mathbb{Z}$, we write $\llbracket a,b \rrbracket := \{a, a+1, \ldots, b\}$ and $\llbracket n \rrbracket := \llbracket 1,n \rrbracket$ with $n \in \mathbb{N}$. In this paper, we shall take $\mathcal{X} = \llbracket n \rrbracket$ unless otherwise specified. 

Let $\ell^2(\pi)$ be the Hilbert space weighted by $\pi$,  with the inner product as 
$$\langle f,g \rangle_\pi := \sum_{x\in \mathcal{X}} f(x)g(x)\pi(x),$$
for $f,g: \mathcal{X} \to \mathbb{R}$. We write $\|f\|^2_\pi = \langle f,f \rangle_\pi$ to be the $\ell^2(\pi)$-norm of $f$. The zero-mean subspace is defined as $\ell^2_0(\pi) := \{f \in \ell^2(\pi): \pi(f) = 0\}$.

Define $\mathcal{L} = \mathcal{L}(\mathcal{X})$ to be the set of all transition matrices on $\mathcal{X}$. For any given $\pi \in \mathcal{P}(\mathcal{X})$, we use $\mathcal{S}(\pi) \subseteq \mathcal{L}$ to denote the set of all $\pi$-stationary transition matrices. For any $P \in \mathcal{S}(\pi),$ it must satisfy $\pi P = \pi$. Similarly, we let $\mathcal{L}(\pi) \subseteq \mathcal{L}$ be the set of all $\pi$-reversible matrices where $P \in \mathcal{L}(\pi)$ implies $\pi(x) P(x,y) = \pi(y) P(y,x)$ holds for all $x,y \in \mathcal{X}$. 

For $P\in \mathcal{S}(\pi)$, $P^*$ is defined to be the time-reversal of $\ell^2(\pi)$-adjoint of $P$. We thus have $P\in \mathcal{L}(\pi)$ if and only if $P^* = P$.

The transition matrices $P \in \mathcal{L}$ can also be viewed as an operator on $\ell^2(\pi)$. Then 
$$Pf(x) = \sum_{y \in \mathcal{X}} P(x,y) f(y)$$
is also a function in $\ell^2(\pi)$. 

For any bounded linear map $T : H_1 \to H_2$ between two Hilbert spaces $H_1, H_2$, we define the operator norm as 
    $$\|T\|_{H_1 \rightarrow H_2} := \sup_{x\neq 0} \frac{\|Tx\|_{H_2}}{\|x\|_{H_1}}.$$
In particular, the operator norm for $P \in \mathcal{L}$ is $\|P\|_{\ell^2(\pi) \rightarrow \ell^2(\pi)}$.

With any $\pi$-reversible $P$ on a finite state space, all eigenvalues are real and lie in $[-1,1]$. We write the distinct eigenvalues in non-increasing order as
$$
1=\lambda_1(P)\ > \lambda_2(P)\ > \dots\ > \lambda_k(P)\ \geq\ -1, \quad 1 \leq k \leq n,
$$
and we denote the set of all distinct eigenvalues of $P$ as 
$$\mathrm{spec}(P) := \{\lambda_1(P), \ldots, \lambda_k(P)\}.$$
Finally, we use $I_k$ to denote the identity matrix of size $k \times k$. If the dimension is clear, we shall drop the subscript and simply use $I$ instead.

\subsection{Group actions}

We now introduce the idea of group actions, which will play a fundamental role in the construction of the proposed samplers. A (left) action of a
group $\mathcal{G}$ on $\mathcal{X}$ is a map $\mathcal{G}\times\mathcal{X}\to\mathcal{X}, (g,x)\mapsto gx,$ such that, for all $g,h\in\mathcal{G}$ and $x\in\mathcal{X}$,
\[ex=x, \qquad (gh)x=g(hx),\]
where $e$ denotes the identity element of $\mathcal{G}$. The action partitions $\mathcal{X}$ into its orbits 
$$\mathcal{O}(x) := \{gx : g \in \mathcal{G}\},$$
and the collection of all orbits is given by $\mathcal{X}/\mathcal G$. The stabiliser of $x$ is then defined by
$$H(x) := \{g \in \mathcal{G}: gx = x\},$$
and for each $y \in \mathcal{O}(x)$,
$$S_y(x) := \{g \in \mathcal{G}: gx = y\}.$$
By the orbit--stabiliser theorem, $S_y(x)$ is a left coset of $H(x)$, so $|S_y(x)|$ is constant across all $y \in \mathcal{O}(x)$.

As outlined in \cite{choi2025groupaveragedmarkovchainsmixing}, we aim to augment $P\in \mathcal{S}(\pi)$ with some suitable group action of $\mathcal{G}$. Formally, at some given state $x\in \mathcal{X}$, we pick $g \in \mathcal G$ with probability 
$$w_x(g) = \frac{\pi(gx)}{\sum_{g\in \mathcal{G}} \pi(gx)},$$ 
and apply the chosen $g$ before applying the sampler $P$ and on the result of $P$.

However, when $|\mathcal{G}|$ is large, direct sampling of the group element $g$ becomes computationally difficult. Hence, instead of working on the group $\mathcal{G}$ itself, we study orbit refreshers defined on the state space $\mathcal{X}$. These are auxillary $\pi$-stationary transition kernels that reshuffle the current state space within its orbit, effectively simulating the effect of sampling $g \in \mathcal{G}$ according to $w_x$ without leaving $\mathcal{X}$. 

We now introduce several such samplers.

\subsection{Gibbs sampler}
Let $G$ denote the orbit refresher kernel on $\mathcal{X}$. By drawing $g \in \mathcal{G}$ with probability $w_x$, we have the formulation of $G$ below.\\

\begin{prop} \label{propG}
For any $x \in \mathcal{X}$, the group--weighted refresher kernel satisfies
\begin{equation} \label{matrixG}
G(x,y) = \begin{cases}
\frac{\pi(y)}{\pi(\mathcal{O}(x))},& \text{for } y \in \mathcal{O}(x),\\
0, & \text{otherwise},
\end{cases}
\end{equation}
where $\pi(\mathcal{O}(x)) := \sum_{z \in \mathcal{O}(x)} \pi(z)$. In particular, the group--based construction coincides with the orbit Gibbs kernel.
\end{prop}

\begin{proof}
Suppose for a given $x \in \mathcal{X}$, we draw $g \in \mathcal{G}$ according to the weights $w_x(g)$.

If $y \notin \mathcal{O}(x),$ then necessarily $G(x,y) = 0$. Otherwise, the transition probability from $x$ to $y \in \mathcal{O}(x)$ is
$$G(x,y) = \sum_{g \in S_y(x)} w_x(g) = \frac{|S_y(x)|\,\pi(y)}{\sum_{z \in \mathcal{O}(x)} |S_z(x)|\,\pi(z)}$$

As $|S_y(x)|$ is independent of $y$, we have the matrix $G$ as
given in \eqref{matrixG}.
\end{proof}

It can then be verified that $G$ is $\pi$-stationary and reversible. It is also an idempotent projection, that is, $G^2 = G$.

\subsection{Metropolis-Hastings and Barker sampler}

Another way to draw $g$ is by running a one–step Metropolis–Hastings move on $\mathcal{G}$.
Given $x$, propose $g$ uniformly from $\mathcal{G}\setminus H(x)$ assuming that we start at the group identity $e$.\\

\begin{prop} \label{propM}
The induced Metropolis-Hastings kernel on $\mathcal{X}$ is
\begin{equation} \label{eqnM}
    M(x,y) = 
    \begin{cases}
    \frac{1}{|\mathcal{O}(x)|-1}\alpha(x,y), &y\in\mathcal{O}(x), y \neq x,\\
    1 - \sum_{y\neq x}M(x,y), & y = x,\\
    0, & \text{otherwise},
    \end{cases}
    \end{equation}
where $\alpha(x,y)=\min\{1,\pi(y)/\pi(x)\}$. If $|\mathcal{O}(x)| = 1$ then $M(x,y) = 0$ for $y\neq x$. This proposal is equivalent to uniformly proposing $y$ within $\mathcal{O}(x)\setminus\{x\}$ and accepting it according the MH rule.
\end{prop}

\begin{proof}
With the proposal 
$$\widetilde Q_x(g)=
\begin{cases}
    \frac{1}{|\mathcal{G}|-|H(x)|}, & g \notin H(x),\\
    0, & \text{otherwise},
\end{cases}$$

and acceptance $\widetilde{\alpha}_x(g)=\min\{1,\pi(gx)/\pi(x)\}$, the one–step MH kernel on $\mathcal{G}$ is

$$\widetilde M_x(g)=\widetilde Q_x(g)\,\widetilde{\alpha}_x(g).$$
To jump from $x$ to a different state $y\neq x$,
$$M(x,y)=\sum_{g\in S_y(x)} \widetilde M_x(g)
=\frac{|S_y(x)|}{|\mathcal{G}|-|H(x)|}\,\min\Big\{1,\frac{\pi(y)}{\pi(x)}\Big\}.$$

By orbit–stabiliser, we have that $|S_y(x)|=|H(x)|$ for all $y\in\mathcal{O}(x)$, and
$$|\mathcal{G}|=|H(x)\|\mathcal{O}(x)|,$$ and hence 
$$M(x,y)=\frac{1}{|\mathcal{O}(x)|-1}\,\min\Big\{1,\frac{\pi(y)}{\pi(x)}\Big\}, \qquad y\in\mathcal{O}(x),\ y\neq x.$$
Setting the diagonals to enforce probability conservation gives the expression given by \eqref{eqnM}. 

If $|\mathcal{O}(x)| = 1$, then $M(x,y) = 0$ for all $y \neq x$ since every group element is in the stabiliser of $x$.
\end{proof}

In fact, the kernel $M$ proposed above is a generalised case of the double-flip move in \cite{ying2022doubleflipisingmodels}. In their construction, the group action is generated by a single involution, coupled with a Metropolis-correction step. This would result in each orbit having size 2, and $M$ being exactly formed by $2\times 2$ blocks. 

A similar kernel can be constructed using the Barker proposal, as defined in \cite{Barker}, with acceptance-rejection ratio 
$$\widetilde{\alpha}^B_x(g) = \frac{\pi(gx)}{\pi(x) + \pi(gx)},$$
which gives us 
\begin{equation} \label{matrixB}
    B(x,y) = 
    \begin{cases}
    \frac{1}{|\mathcal{O}(x)|-1}\alpha^B(x,y), &y\in\mathcal{O}(x), y \neq x,\\
    1 - \sum_{y\neq x}B(x,y), & y = x,\\
    0, & \text{otherwise},
    \end{cases}
\end{equation}
where $\alpha^B(x,y) = \frac{\pi(y)}{\pi(x) + \pi(y)}$. Again if $|\mathcal{O}(x)| = 1$, $B(x,y) = 0$ for all $y \neq x$. Note that if $\mathcal{G}$ admits orbits all of size at most 2, then $B = G$. 

In fact, in almost all cases, $G$ is the limit of $B^k$ and $M^k$:\\
\begin{prop} \label{limMB}
    Assume that the same group $\mathcal{G}$ and its associated group action is used in defining $B$, $M$ and $G$.
    If $M$ does not have a deterministic 2-cycle on any of its orbits, then
    $$\lim_{i\rightarrow \infty} M^i = G.$$
    The same limit holds for any $B$, that is,
    $$\lim_{i\rightarrow \infty} B^i = G.$$
\end{prop}

\begin{proof}
    For all three kernels, they can be written in block diagonal form in terms of their orbits. For example, if $|\mathcal{X}/\mathcal{G}| = k$, we have that
    $$G = \text{diag}(G_{\mathcal{O}_1} , \dots, G_{\mathcal{O}_k})$$
    and each $G_{\mathcal{O}_k}$ has identical rows 
    $$\pi_{\mathcal{O}_k} = \frac{1}{\pi(\mathcal{O}_k)}\bigg(\pi(x_1), \dots, \pi(x_k)\bigg).$$

    If we similarly decompose $B$, one may then verify that $\pi_{\mathcal{O}_k} B_{\mathcal{O}_k} = \pi_{\mathcal{O}_k}$, and that each block $B_{\mathcal{O}_k}$ is always ergodic. Hence, $\lim_{i \rightarrow \infty}B_{\mathcal{O}_k}^i = G_{\mathcal{O}_k}$, and naturally, 
    $$\lim_{i\rightarrow \infty} B^i = \lim_{i\rightarrow \infty} \text{diag}(B_{\mathcal{O}_1}^i , \dots, B_{\mathcal{O}_k}^i) = \text{diag}(G_{\mathcal{O}_1} , \dots, G_{\mathcal{O}_k}) = G.$$

    It also holds true for $M$, that $\pi_{\mathcal{O}_k} M_{\mathcal{O}_k} = \pi_{\mathcal{O}_k}$. However, $M_{\mathcal{O}_k}$ is ergodic if and only if 
    $$M_{\mathcal{O}_k} \neq \begin{bmatrix}
    0 & 1\\
    1 & 0
    \end{bmatrix}.$$
    So in the case where the deterministic 2-cycle does not occur, a similar conclusion to that of $B$ will be true, and $\lim_{i\rightarrow \infty} M^i = G$.
\end{proof}

\begin{remark}
Although our construction of group averaging is motivated by group
actions, many of the results below depend only on the induced orbit
partition. Indeed, every partition of a finite state space can be realised
as the orbit partition of a finite group action: given a partition
$(\mathcal{O}_i)_{i=1}^k$, one may take a permutation whose disjoint cycles
are precisely the non-singleton blocks $\mathcal{O}_i$, with singleton
blocks fixed, and consider the cyclic group generated by this permutation.
We retain the group-action formulation because such partitions arise
naturally from symmetries and because group actions provide a convenient
mechanism for constructing and sampling within the orbits. Whenever a
result depends only on the induced partition rather than on the underlying
group structure, we make this dependence explicit.
\end{remark}

\subsection{An intuitive orbit perspective to understand the improvement in mixing of $GPG,MPM,BPB$ over $P$}

The crux of the group-averaged approach is to augment the original $P$ with a group $\mathcal{G}$ that acts on $\mathcal{X}$. In doing so, this induces a partition of $\mathcal{X}$ based on the group orbits. Precisely, suppose that $|\mathcal{X}/\mathcal{G}| = k$, and so one can write that
$$\mathcal{X} = \bigcup_{i=1}^k\mathcal{O}_i.$$

We note that $G,M,B$ facilitates within orbit transitions (e.g. from $\mathcal{O}_i$ to $\mathcal{O}_i$), which might be hard to reach using $P$ only. On the other hand, the original $P$ is capable to facilitates both within orbit and cross orbit transitions (e.g. from $\mathcal{O}_i$ to $\mathcal{O}_j$ with $i \neq j$). Thus, using $GPG, MPM, BPB$ are capable to enhance within orbit transitions over the original $P$.

We now recall three important notions from \citet{Jerrum_2004}. The first one is the notion of projection chain $\overline{P}$ induced by the partition $(\mathcal{O}_i)_{i=1}^k$, where $\overline{P}: \llbracket k \rrbracket \times \llbracket k \rrbracket \to [0,1]$ is defined to be
\begin{align} \label{Pbar}
    \overline{P}(i,j) := \frac{1}{\pi(\mathcal{O}_i)} \sum_{\substack{x \in \mathcal{O}_i \\ y \in \mathcal{O}_j}} \pi(x) P(x,y),
\end{align}
with stationary distribution $\overline{\pi} = (\pi(\mathcal{O}_1), \dots, \pi(\mathcal{O}_k))$. Note that the dependency of $\overline{P}$ on $\mathcal{G}$ is suppressed.

The second one is the notion of restriction chains $P_1,P_2,\ldots,P_k$ induced by the partition $(\mathcal{O}_i)_{i=1}^k$, with $P_i: \mathcal{O}_i \times \mathcal{O}_i \to [0,1]$ defined by 
\begin{align} \label{P_i}
    P_i(x,y) := \begin{cases}
        P(x,y), & \text{if } x\neq y,\\
        1 - \displaystyle\sum_{z\in \mathcal{O}_i\setminus \{x\}} P(x,z), & \text{if } x = y,
    \end{cases}
\end{align}
and stationary distribution $\pi_i(x) = \pi(x)/\pi(\mathcal{O}_i)$. Note that the dependency of $P_i$ on $\mathcal{G}$ is suppressed.

The third one is the notion of $\gamma(P)$ induced by the partition $(\mathcal{O}_i)_{i=1}^k$:
\begin{align} \label{gamma(P)}
    \gamma(P) := \max_{i \in \llbracket k \rrbracket}~\max_{x \in \mathcal{O}_i} \sum_{y \in \mathcal{X}\setminus \mathcal{O}_i} P(x,y).
\end{align}

Again, the dependency of $\gamma$ on $\mathcal{G}$ is suppressed. Analogously, we write $\overline{GPG}, ((GPG)_i)_{i=1}^k$ and $\gamma(GPG)$ to be the projection chain, restriction chains and $\gamma$ respectively of $GPG$ induced by the partition $(\mathcal{O}_i)_{i=1}^k$. We denote similar objects for $BPB,MPM$ as well. 

Below, we attempt to compare the samplers $P$ and $GPG$ and their above-mentioned counterparts in terms of right spectral gap.\\

\begin{prop} \label{projchain_equal}
    The projection chain of $P$ and $GPG$ induced by the partitions $(\mathcal{O}_i)_{i=1}^k$ are identical, or equivalently 
    $$\overline{GPG} = \overline{P}.$$
\end{prop}

\begin{proof}
    For any $x, y \in \mathcal{X}$, we have that
    \begin{equation} \label{GPGmatrix}
        GPG(x,y) = \frac{\pi(y)}{\pi(\mathcal{O}(x))\pi(\mathcal{O}(y))} \sum_{\substack{z \in \mathcal{O}(x)\\ w \in \mathcal{O}(y)}} \pi(z) P(z,w)
    \end{equation}
    With that, 
    \begin{align*}
        \overline{GPG}(i,j) &= \frac{1}{\pi(\mathcal{O}_i)} \sum_{\substack{x \in \mathcal{O}_i\\ y \in \mathcal{O}_j}}\bigg[ \pi(x)\cdot GPG(x,y)\bigg]\\ 
        &= \frac{1}{\pi(\mathcal{O}_i)} \sum_{\substack{x,z \in \mathcal{O}_i\\ y,w \in \mathcal{O}_j}}\bigg[ \frac{\pi(x)\cdot \pi(y)\cdot \pi(z)}{\pi(\mathcal{O}_i)\pi(\mathcal{O}_j)} P(z,w) \bigg]\\
        &= \frac{1}{\pi(\mathcal{O}_i)} \sum_{\substack{z \in \mathcal{O}_i\\ w \in \mathcal{O}_j}}\bigg[ \pi(z)\cdot P(z,w)\bigg] = \overline{P}(i,j)
    \end{align*}
\end{proof}

\begin{prop} \label{gammaequal}
    The inequality 
    $$\gamma(P) \geq \gamma(GPG)$$
    holds true for any choice of $\pi$-stationary $P$. Equality holds when the maximum-achieving orbit $\mathcal{O}_i$ satisfies the property that $\sum_{y \notin \mathcal{O}_i} P(x,y)$ is equal for all $x \in \mathcal{O}_i$.
\end{prop}

\begin{proof}
    Fix orbit $\mathcal{O}_i$ and take $x \in \mathcal{O}_i$. Then
    \begin{align*}
        \sum_{y \notin \mathcal{O}_i} GPG(x,y) &= \sum_{y \notin \mathcal{O}_i}\Bigg[ \frac{\pi(y)}{\pi(\mathcal{O}_i)\pi(\mathcal{O}(y))} \sum_{\substack{z\in \mathcal{O}_i \\ w \in \mathcal{O}(y)}} \pi(z)\cdot P(z,w)\Bigg]\\
        &= \sum_{j \neq i} \sum_{y \in \mathcal{O}_j} \Bigg[ \frac{\pi(y)}{\pi(\mathcal{O}_i)\pi(\mathcal{O}_j)} \sum_{\substack{z\in \mathcal{O}_i \\ w \in \mathcal{O}_j}} \pi(z)\cdot P(z,w)\Bigg]\\
        &= \sum_{j \neq i}\Bigg[ \frac{1}{\pi(\mathcal{O}_i)} \sum_{\substack{z\in \mathcal{O}_i \\ w \in \mathcal{O}_j}} \pi(z) P(z,w)\Bigg] \\
        &= \sum_{z\in \mathcal{O}_i}\bigg[ \frac{\pi(z)}{\pi(\mathcal{O}_i)} \sum_{w \notin \mathcal{O}_i} P(z,w)\bigg],
    \end{align*}
    is in fact, independent of $x$. Furthermore, since it is a convex combination, it must be that
    $$ \sum_{y \notin \mathcal{O}_i} GPG(x,y) \leq \max_{x\in \mathcal{O}_i} \sum_{y \notin \mathcal{O}_i} P(x,y)$$
    for any choice of $x$. Taking maximum again over all possible orbits yields the desired inequality. 

    The same inequality also shows that for equality to hold, we require the inner sum to be constant across all $z\in \mathcal{O}_i$ on the maximum-achieving orbit.
\end{proof}

\begin{prop} \label{eigenGPG_i}
    For each $i \in \llbracket k\rrbracket$, the restriction chain $(GPG)_i$ defined on $\mathcal{O}_i$ has eigenvalues $\lambda_1 = 1,\ \lambda_2 = 1-\overline{a}_i$, where
    $$a_i(x) := \sum_{y \in \mathcal{O}_i} P(x,y)\quad \text{and} \quad \overline{a}_i := \sum_{x\in \mathcal{O}_i} \frac{\pi(x)}{\pi(\mathcal{O}_i)} \sum_{y \in \mathcal{O}_i} P(x,y) = \mathbb{E}_{\pi_i} a_i(x)$$
\end{prop}

\begin{proof}
    On the orbit $\mathcal{O}_i$, 
    $$(GPG)_i(x,y) = \begin{cases}
        \displaystyle\frac{\pi(y)}{\pi(\mathcal{O}_i)} \overline{a}_i, & \text{for } x \neq y,\\
        1 - \displaystyle\sum_{y \neq x} \frac{\pi(y)}{\pi(\mathcal{O}_i)} \overline{a}_i, & \text{for } x = y. 
    \end{cases}$$
    In fact, $(GPG)_i = (1-\overline{a}_i)I + \overline{a}_i G_i$, and given that the rank of $G_i = 1$, the eigenvalues must be  
    $$\lambda_1 = 1,\ \lambda_2 = 1-\overline{a}_i.$$
\end{proof}

It is natural to ask whether the improvement guaranteed by group averaging extends to the restriction chains $((GPG)_i)_{i=1}^k$ individually. The following example shows that this need not hold in general.

Consider $\pi = (0.3, 0.3, 0.4)$, $\mathcal{G} = \{e, (1,2)\}$ and 
$$P(x,y) = \begin{pmatrix}
    0 & 0.4 & 0.6\\
    0.4 & 0 & 0.6\\
    0.45 & 0.45 & 0.10
\end{pmatrix}.$$
Then on the orbit $\{1,2\}$,
$$
P_1 = \begin{pmatrix}
    0.6 & 0.4\\
    0.4 & 0.6
\end{pmatrix}
$$
and 
$$
GPG = \begin{pmatrix}
    0.2 & 0.2 & 0.6\\
    0.2 & 0.2 & 0.6\\
    0.45 & 0.45 & 0.10
\end{pmatrix}, \quad
(GPG)_1 = \begin{pmatrix}
    0.8 & 0.2\\
    0.2 & 0.8
\end{pmatrix}.
$$
The corresponding eigenvalues are $\lambda_2(P_1)=0.2$ and $\lambda_2((GPG)_1)=0.6$, so the local spectral gap of $(GPG)_1$ is strictly smaller than that of $P_1$.

Hence, while $GPG$ globally improves or preserves the overall spectral properties of $P$, the same is not necessarily true within each individual orbit.

\subsection{Additive group-averaged Markov chains}

In previous sections, we have introduced multiplicative group-averaged Markov chains $GPG, MPM, BPB$, which are respectively based on the Gibbs kernel, Metropolis-Hastings kernel and Barker kernel.

Another class of group-averaged Markov chains is additive group-averaged Markov chains, which are defined to be
\begin{align*}
    A(G,P) &:= \dfrac{1}{2}(G + P), \\
    A(M,P) &:= \dfrac{1}{2}(M + P), \\
    A(B,P) &:= \dfrac{1}{2}(B + P),
\end{align*}
that we call respectively the additive group-averaged Gibbs sampler, additive group-averaged Metropolis-Hastings sampler and additive group-averaged Barker sampler. More generally, the additive Gibbs mixture $\alpha P+(1-\alpha)G,$ for $\alpha\in[0,1]$ is studied in detail by \cite{additive}, where it is interpreted through a lifted Markov chain and the choice of orbit partition is analysed under the Frobenius norm and KL divergence.

In the special case where $\mathcal{G} = \{e\}$, we see that $G = M = B = I$, and hence $A(G,P) = A(M,P) = A(B,P) = (1/2)(I+P)$, the lazified version of $P$.

Additive mixtures of this form also appear in existing symmetry-based samplers. For example, the double-flip Swendsen–Wang algorithm of \cite{ying2022doubleflipisingmodels} is a special case, with $P$ given by the usual Swendsen-Wang dynamics and $M$ given by their Metropolis double-flip move. 

We first list several properties that are related to these additive samplers.\\

\begin{prop} \label{Gprop}
    For the group-orbit samplers $G$ (and also, $B$ and $M$), they satisfy the following properties:
    \begin{itemize}
        \item $\overline{G} = I_k$.
        \item $G_i$ is identical to the block diagonal matrix of $G$ on the orbit $\mathcal{O}_i$.
        \item $\gamma(G) = 0.$\\
    \end{itemize}
\end{prop}

\begin{prop} \label{convexprop}
    Consider two $\pi$-stationary samplers $P$ and $Q$, and the same partition of orbits $(\mathcal{O}_i)_{i=1}^k$. Then for any $\alpha \in [0,1]$,
    $$\overline{\alpha P + (1-\alpha) Q} = \alpha \overline{P} + (1-\alpha) \overline{Q}$$
    and 
    $$(\alpha P + (1-\alpha)Q)_i = \alpha P_i + (1-\alpha) Q_i.$$
    Furthermore, if $Q$ is one of the orbit samplers $G, M$ or $B$, 
    $$\gamma(\alpha P + (1-\alpha) Q) = \alpha \gamma(P).$$
\end{prop}

\begin{corollary} \label{mixturedecomp}
    Let $Q$ be any of the orbit samples $G, M$ or $B$. Then the mixture $K_\alpha(Q) = \alpha P + (1-\alpha)Q$ has the following properties:
    \begin{itemize}
        \item $\lambda_2(\overline{K_\alpha(Q)}) = 1 - \alpha + \alpha\lambda_2(\overline{P}).$
        \item $\lambda_2\big([K_\alpha(Q)]_i\big) \leq \alpha\lambda_2(P_i) + (1-\alpha)\lambda_2(Q_i).$ In particular if $Q = G$, then $\lambda_2\big([K_\alpha(Q)]_i\big) = \alpha\lambda_2(P_i)$.
        \item $\gamma(K_\alpha(Q)) = \alpha\gamma(P)$
    \end{itemize}
\end{corollary}

Despite these explicit properties of the projection and restriction chains, the right spectral gap of $K_\alpha(Q)$ cannot in general be ordered uniformly relative to that of $P$.

\begin{prop}
    Let $\lambda(P) := 1 - \lambda_2(P)$ be the right spectral gap of a $\pi$-reversible sampler $P$. In general, there exists no uniform ordering of $\lambda(K_\alpha(Q))$ and $\lambda(P)$. 
\end{prop}

\begin{proof}
    Let $\mathcal{X} = \llbracket 4 \rrbracket$. Set $\pi$ to be the uniform distribution on $\mathcal{X}$ and $P = 1/2(I + \Pi).$ One can verify that $\lambda(P) = 1/2.$
    
    \textbf{Example where $Q$ improves the gap:}
    Consider the trivial partition where the entire space $\mathcal{X}$ is a single orbit. Then $G = \Pi$,
    $$M = \begin{pmatrix}
        0 & 1/3 & 1/3 & 1/3\\
        1/3 & 0 & 1/3 & 1/3\\
        1/3 & 1/3 & 0 & 1/3\\
        1/3 & 1/3 & 1/3 & 0
    \end{pmatrix}, \quad
    B = \begin{pmatrix}
        1/2 & 1/6 & 1/6 & 1/6\\
        1/6 & 1/2 & 1/6 & 1/6\\
        1/6 & 1/6 & 1/2 & 1/6\\
        1/6 & 1/6 & 1/6 & 1/2
    \end{pmatrix}.
    $$

    One can verify that that the right spectral gaps of $K_{1/2}(Q)$ are 3/4, 11/12 and 7/12 for $Q = G,M,B$ respectively, all of which are strictly larger than $\lambda(P) = 1/2.$

    \textbf{Example where $Q$ worsens the gap:}
    Take the orbits $\{1,2\}, \{3,4\}$. Then 
    $$G = B = \begin{pmatrix}
        1/2 & 1/2 & 0 & 0\\
        1/2 & 1/2 & 0 & 0\\
        0 & 0 & 1/2 & 1/2\\
        0 & 0 & 1/2 & 1/2
    \end{pmatrix},\ 
    M = \begin{pmatrix}
        0 & 1 & 0 & 0\\
        1 & 0 & 0 & 0\\
        0 & 0 & 0 & 1\\
        0 & 0 & 1 & 0
    \end{pmatrix}.
    $$
    One can verify that the right spectral gaps of $K_{1/2}(Q)$ for $Q = G,M,B$ are all 1/4, which is lesser than $\lambda(P) = 1/2$.
\end{proof}

\section{Comparison of absolute spectral gap} \label{section:specgap}
In this section, we assume $P\in \mathcal{L}(\pi)$ to be an ergodic time-reversible sampler of $\pi$. We attempt to compare the absolute spectral gap of $P$, $GPG$, $BPB$ and $MPM$. Note that the multiplicative samplers are $\pi$-stationary, and admit real spectra that lies within $[-1, 1]$.

For $T \in \mathcal{L}(\pi)$, we write 
$$\text{Fix}(T) = \{ f \in \ell^2(\pi): Tf = f\},$$
and Fix$(T)^\perp$ to be the orthocomplement. 

We then define the second-largest eigenvalue in modulus (SLEM) to be
$$\rho(T) := \max\{|\lambda_{2}(P)|,\ |\lambda_k(P)|\}, \quad (\lambda_i(T))_{i=1}^k \in \mathrm{spec}(T),$$
and 
$$\lambda(P)\ :=\ 1-\rho(P)$$
to be the absolute spectral gap in this section.

For any self-adjoint $T$, we also have that $\rho(T) = \|T|_{\text{Fix}(T)^\perp}\|_{\ell^2(\pi) \rightarrow \ell^2(\pi)}$.

An equivalent definition of SLEM given in the Rayleigh-Ritz form is 
\begin{equation} \label{SLEMdef}
    \rho(T) = \sup_{f \neq 0, f \in \ell^2_0(\pi)}\frac{|\langle f, Tf\rangle_\pi|}{\langle f,f \rangle_\pi}.
\end{equation}

Hence, comparisons involving spectral gap can be equivalently computed by the comparisons in $\lambda_2$ and $\lambda_k$ of the different samplers.

\subsection{Comparing original and group-averaging kernels}
We first compare $P$ to the Gibbs-orbit sampler $GPG$.\\

\begin{lemma} \label{Gortho}
    The Gibbs orbit kernel $G$  is an orthogonal projection onto the subspace 
    $$S := \{f \in \ell^2(\pi) : f(x) = f(y) \text{ for } y \in \mathcal{O}(x)\}.$$
    In other words, the subspace $S$ contains only functions that are constant on each orbit of $\mathcal{G}$.
\end{lemma}

\begin{proof}
    Take any $f \in \ell^2(\pi)$, and for any $x \in \mathcal{X}$,
    $$Gf(x) = \sum_{y\in \mathcal{X}} f(y)G(x,y) = \frac{1}{\pi(\mathcal{O}(x))}\sum_{y\in \mathcal{O}(x)} f(y)\pi(y).$$
    The same result holds for any $x' \in \mathcal{O}(x)$, hence $Gf$ must be constant on orbits. 
\end{proof}

\begin{prop} \label{SLEMGPG}
    Let $P$ be self-adjoint and $G$ be the Gibbs kernel defined in Proposition \ref{propG}. Then we always have $\rho(GPG) \leq \rho(P)$, or equivalently, $\lambda(GPG) \geq \lambda(P)$.

    Moreover, equality holds if and only if a SLEM-achieving eigenfunction lies in $S$, the subspace projected by $G$ defined in Lemma \ref{Gortho}.
\end{prop}

\begin{proof}
    By the Rayleigh-Ritz characterization on $\ell^2_0(\pi)$ in \eqref{Pbar},
    
    $$\rho(GPG) = \sup_{f\neq 0}\frac{|\langle f,\,GPG f\rangle_\pi|}{\langle f,f\rangle_\pi}
    = \sup_{f\neq 0}\frac{|\langle Gf,\,P\,Gf\rangle_\pi|}{\langle f,f\rangle_\pi},$$
    using $G=G^*$. Since $G$ is a contraction on $\ell^2(\pi)$, $\|Gf\|_\pi\le \|f\|_\pi$, so
    
    $$\frac{|\langle Gf,\,P\,Gf\rangle_\pi|}{\langle f,f\rangle_\pi}
    \leq \frac{|\langle Gf,\,P\,Gf\rangle_\pi|}{\langle Gf,Gf\rangle_\pi}.$$
    
    Writing $u=Gf$, we have $u\in S\cap \ell^2_0(\pi)$ and thus
    $$
    \rho(GPG)\ \leq \sup_{\substack{u\in S\cap \ell^2_0(\pi)\\ u\neq 0}}
    \frac{|\langle u,Pu\rangle_\pi|}{\langle u,u\rangle_\pi} \leq 
    \sup_{\substack{u\in \ell^2_0(\pi)\\ u\neq 0}}
    \frac{|\langle u,Pu\rangle_\pi|}{\langle u,u\rangle_\pi}
    = \rho(P),
    $$
    where the last inequality is a result of taking supremum over a larger set.
    
    For equality, note that the inequalities are tight if and only if there exists $u^\star\in S\cap \ell^2_0(\pi)$ attaining the $P$-supremum. That is equivalent to having a SLEM-achieving eigenfunction of $P$ belonging to $S$.
\end{proof}

Similar results hold for the samplers $BPB$ and $MPM$ as well. \\

\begin{prop} \label{SLEMMPM}
    Let $P$ be self-adjoint and $M$, $B$ be the MH-orbit and Barker-orbit defined in \eqref{eqnM} and \eqref{matrixB}. Then $\rho(MPM) \leq \rho(P)$ and $\rho(BPB) \leq \rho(P)$. Equivalently, the absolute spectral gap of both $MPM$ and $BPB$ are no worse than that of $P$.
\end{prop}

\begin{proof}
    Again by the Rayleigh-Ritz definition of SLEM in \eqref{SLEMdef}, for $f \in \ell^2_0(\pi)$,
    $$\rho(MPM) =\sup_{f\neq 0}\frac{|\langle f,MPM f\rangle_\pi|}{\langle f,f\rangle_\pi} =\sup_{f\neq 0}\frac{|\langle Mf, P\,Mf\rangle_\pi|}{\langle f,f\rangle_\pi},$$
    
    where the last equality is a result of $M$ being self-adjoint on $\ell^2(\pi)$. Additionally, $M$ is a contraction, that is, $\|Mf\|_\pi\le \|f\|_\pi$, and hence
    $$
    \frac{|\langle Mf, P\,Mf\rangle_\pi|}{\langle f,f\rangle_\pi}
    \;\le\;
    \frac{|\langle Mf, P\,Mf\rangle_\pi|}{\langle Mf,Mf\rangle_\pi}
    \;\le\;
    \sup_{u\neq 0}\frac{|\langle u,Pu\rangle_\pi|}{\langle u,u\rangle_\pi}
    =\rho(P).
    $$
    Taking the supremum over $f$ yields $\rho(MPM)\le \rho(P)$.
    The same argument holds when we replace $M$ with $B$.
\end{proof}

The proof of Proposition \ref{SLEMMPM} also reveals that equality holds if and only if there exists an eigenfunction $f \in \ell^2_0(\pi)$ achieving the SLEM of $P$ such that
$$\|Mf\|_\pi = \|f\|_\pi \quad \text{and}\quad \frac{|\langle Mf, PMf\rangle_\pi|}{\langle Mf, Mf \rangle_\pi} = \frac{|\langle f, Pf \rangle_\pi|}{\langle f,f \rangle_\pi}.$$
Equivalently, equality holds when the SLEM–achieving eigenfunction $f$ of $P$ is also an eigenfunction of $M$ with eigenvalue $\pm 1$.

\paragraph{Case of $M$ with eigenvalue $+1$.}
If $M$ does not admit the eigenvalue $-1$, equality requires $Mf = f$.  
This condition is satisfied when $f$ is constant on each orbit of the group action, but not globally constant.  
In particular, when the entire state space forms a single orbit, the only function satisfying both $Mf = f$ and $f \in \ell^2_0(\pi)$ is the zero function, so equality cannot occur in this case.

\paragraph{Case of $M$ with eigenvalue $-1$.}
The situation $Mf = -f$ arises only under a two–cycle, where $M$ acts as a deterministic flip between two states $x_1, x_2$ with equal stationary weights, i.e. $\pi(x_1) = \pi(x_2)$.  
In this case, an antisymmetric eigenfunction supported on that orbit (e.g. $f(x_1) = 1, f(x_2) = -1$) yields equality.

\paragraph{Barker proposal.}
For the Barker kernel $B$, the acceptance probability $\alpha^B$ ensures the presence of self–loops and hence cannot have eigenvalue $-1$.  
Therefore, equality in $\rho(BPB) = \rho(P)$ can occur only when the SLEM–eigenfunction of $P$ is constant on each orbit but not globally constant.

\subsection{Comparison between different group-averaging kernels}
We now show that under certain conditions, the spectral gap of $GPG$ is never worse than that of $MPM$ or $BPB$. Here, we assume that the same group $\mathcal{G}$, and the same ergodic kernel $P$ is used, with the only difference being the choice of sampler for the group action.\\

\begin{lemma} \label{GM=MG}
    For the same group $\mathcal{G}$, we have that $GM = MG = BG = GB = G$.
\end{lemma}

\begin{proof}
    We first show the equality $MG = G$. For any $x \in \mathcal{X}$ and $z \in \mathcal{O}(x)$,
    $$
    MG(x,z) = \sum_{y\in \mathcal{O}(x)} M(x,y)G(y,z) = \sum_{y\in \mathcal{O}(x)} \frac{\pi(z)}{\pi(\mathcal{O}(y))} M(x,y) = G(x,z).
    $$
    If $z \notin \mathcal{O}(x)$, then equality holds trivially for $MG(x,z) = G(x,z) = 0$.

    For $GM = G$, we show that 
    $$\text{Fix}(M) := \{f : Mf = f\}$$
    is equivalent to $S = \text{Fix}(G)$.
    
    If $f \in S$, then for any $x \in \mathcal{X}$,
    $$
    (Mf)(x) = \sum_{y\in\mathcal{O}(x)} M(x,y)f(y) = f(x) \sum_y M(x,y) = f(x)
    $$
    since $f$ is constant on $\mathcal{O}(x)$. Hence $S \subseteq \text{Fix}(M)$.

    Now for $f \in \text{Fix}(M)$, take any orbit $\mathcal{O}$ from $\mathcal{X}/\mathcal{G}$. Suppose within the orbit $f$ reaches a maximum at $x \in \mathcal{O}$. Then
    $$f(x) = (Mf)(x) = \sum_{y \in \mathcal{O}} M(x,y)f(y),$$
    and so $f(x)$ is some convex combination of $f(y)$ for $y \in \mathcal{O}$. For equality to hold, we must have $f$ being constant on the entire orbit $\mathcal{O}$. Equivalently, $f \in S,$ and so Fix$(M) \subseteq S$. 

    Finally, it must be the case that $MG = G$ since $Gf \in S$, so $M(Gf) = Gf$ for all $f \in \ell^2(\pi)$. The same argument holds by replacing $M$ with $B$. 
\end{proof}

Now let $R = M-G$ be the additive decomposition of $M$. One can verify that $R$ is self-adjoint with the following properties 
\begin{itemize}
    \item $R$ annihilates the subspace $S$; for $f \in S$,
    $$Rf = (M-G)f = 0.$$

    \item $R$ maps onto $S^\perp$, the subspace orthogonal to $S$. 

    \item It acts like $M$ on $S^\perp$, that is $Mf = Rf$ for any $f \in S^\perp$.

    \item The spectrum of $R$ is exactly the spectrum of $M|_{S^\perp} \cup \{0\}$. Equivalently, 
    $$\|R\|_{\ell^2(\pi) \rightarrow \ell^2(\pi)} = \max\{|\lambda_i|: \lambda_i \in \mathrm{spec}(M|_{S^\perp})\} := \theta \leq 1,$$
    with strict inequality if each orbit chain is aperiodic.\\ 
\end{itemize}

\begin{prop} \label{SLEMdiffGPG/MPM}
    We have the inequality 
    $$0 \leq |\lambda_2(MPM)| - |\lambda_2(GPG)| \leq \rho(P|_{S^\perp})(2\theta + \theta^2)$$
    for $G$ and $M$ sharing the same group action $\mathcal{G}$.
\end{prop}

\begin{proof}
    First, for any $f \in S$, $(MPM-GPG)f = 0$, so we restrict our attention only to $S^\perp$. By the triangle inequality,
    $$\|MPM - GPG\|_{\ell^2(\pi) \rightarrow \ell^2(\pi)} \leq \|GPR\|_{\ell^2(\pi) \rightarrow \ell^2(\pi)} + \|RPG\|_{\ell^2(\pi) \rightarrow \ell^2(\pi)} + \|RPR\|_{\ell^2(\pi) \rightarrow \ell^2(\pi)}.$$

    Since $G$ is idempotent, $\|G\|_{\ell^2(\pi) \rightarrow \ell^2(\pi)} = 1$. Then by the submultiplicativity properties of norm, and that $P$ is assumed to be ergodic,
    $$\|MPM - GPG\|_{\ell^2(\pi) \rightarrow \ell^2(\pi)} \leq \rho(P|_{S^\perp})(2\theta + \theta^2).$$
    Using Weyl's inequality, we have that 
    \begin{equation} \label{weyl}
        |\lambda_i(MPM) - \lambda_i(GPG)| \leq \|MPM-GPG\|_{\ell^2(\pi) \rightarrow \ell^2(\pi)} \leq \rho(P|_{S^\perp})(2\theta + \theta^2).
    \end{equation}

    Further, using the fact that Fix($M$) = $S$,
    \begin{align*}
    |\lambda_2(MPM)| &= \sup_{f\in \ell^2_0(\pi)\setminus\{0\}} \frac{|\langle f, MPMf\rangle_\pi|}{\langle f,f \rangle_\pi}\\
    &\geq \sup_{u\in S_0\setminus\{0\}} \frac{|\langle u, MPMu\rangle_\pi|}{\langle u,u \rangle_\pi}\\
    &= \sup_{u\in S_0\setminus\{0\}} \frac{|\langle u, GPGu\rangle_\pi|}{\langle u,u \rangle_\pi} = |\lambda_2(GPG)|
    \end{align*}
    where $S_0 = \ell^2_0(\pi) \cap S$. So 
    $$0 \leq |\lambda_2(MPM)| - |\lambda_2(GPG)| \leq \rho(P|_{S^\perp})(2\theta + \theta^2).$$
\end{proof}

\begin{corollary} \label{SLEMdiffGPG/BPB}
    The inequality $$0 \leq |\lambda_2(BPB)| - |\lambda_2(GPG)| \leq  \rho(P|_{S^\perp})(2\theta + \theta^2)$$ holds with the same argument as in Proposition \ref{SLEMdiffGPG/MPM}, by interchanging $M$ and $B$.
\end{corollary}

While Propositions \ref{SLEMdiffGPG/MPM} and \ref{SLEMdiffGPG/BPB} give a clear ordering of $|\lambda_2|$ for $MPM/BPB$ with $GPG$, such an ordering does not extend to the smallest eigenvalue $\lambda_n$. In particular, when the negative eigenvalues dominate, we do not have strict ordering of SLEM.

A sufficient condition for $\rho(GPG) \leq \rho(MPM)$ is that $P$ has non-negative spectrum. Then 
$$\rho(GPG) = |\lambda_2(GPG)| \leq |\lambda_2(MPM)| = \rho(MPM).$$

However, in cases where $|\mathcal{G}|$ is large, the calculation of $G$ becomes computationally infeasible. Since $M$ and $B$ are much easier to simulate, and both limits tend towards $G$ under most circumstances, we now try to quantify the rate of convergence as we take increasing powers of $M$ or $B$ to approximate $G$.

We hence have the following result:
\begin{prop} \label{SLEMM^kPM^k}
    For all ordered pairs of eigenvalues
    $$|\lambda_i(M^k P M^k) - \lambda_i(GPG)| \leq \rho(P|_{S^\perp})(2\theta^k + \theta^{2k}).$$
    In particular, $0 \leq \lambda_2(M^k P M^k) - \lambda_2 (GPG) \leq \rho(P|_{S^\perp})(2\theta^k + \theta^{2k}) \longrightarrow 0$ if $\theta = \|R\|_{\ell^2(\pi) \rightarrow \ell^2(\pi)} < 1$.
\end{prop}

\begin{proof}
    Observe that 
    $$M^k P M^k - GPG = (G+R^k) P (G+R^k) - GPG = GPR^k + R^kPG + R^kPR^k.$$
    Applying the operator norm and using the properties of subadditivity and submultiplicity, we then have $$\|M^k P M^k - GPG\|_{\ell^2(\pi) \rightarrow \ell^2(\pi)} \leq \rho(P)(2\theta^k + \theta^{2k}).$$
    Weyl's inequality then gives us the desired inequality for the absolute difference in eigenvalues. 

    The inequality $\lambda_2(M^k P M^k) \leq \lambda_2(GPG)$ can then be shown by a similar argument as in Proposition \ref{SLEMdiffGPG/BPB}, noting that $M^kG = G = GM^k$. 
\end{proof}

The result above shows that we observe convergence with a factor of $\theta$, so long as $\theta < 1$. Since $\theta$ is the spectrum of $R = M|_{S^\perp}$, it can never achieve the eigenvalue $1$ since all such eigenfunctions lie in $S$. Again, the eigenvalue $-1$ can only be achieved in the sole case where there exists a degenerate 2-cycle in one of our orbits.  

The above results naturally also extends to that of $B$ as well. However, with $B$, $\theta < 1$ always holds since it cannot have a 2-cycle by design. 

In the Metropolis-Hastings case, we can further characterise $\theta$ in the following manner:

\begin{prop} \label{thetadef}
    Suppose the group action of $\mathcal{G}$ admits $k$ orbits, $\mathcal{O}_1, \dots, \mathcal{O}_k$, not all of size 1. On each orbit, label the elements of $\mathcal{O}_i$ by non-increasing order in terms of $\pi$:

    $$\pi(x_1^{(i)}) \geq \pi(x_2^{(i)}) \geq \dots \geq \pi(x_{m_i}^{(i)}),$$
    where $m_i = |\mathcal{O}_i|$ and $x_k^{(i)} \in \mathcal{O}_i$. Then 
    $$\theta = \rho(M|_{S^\perp}) = \max_{i \in \llbracket k \rrbracket;~ m_i > 1} \bigg\{\bigg| 1-\frac{\pi(\mathcal{O}_i)}{(m_i-1)\pi(x_1^{(i)})} \bigg|,\ \frac{\pi(x_{m_i}^{(i)})}{\pi(x_{m_i-1}^{(i)})(m_i -1)} \bigg\}$$
    where the maximum is taken across all orbits $\mathcal{O}_i$ with $m_i >1$.
\end{prop}

\begin{proof}
    We first look at a single orbit $\mathcal{O}_i$ with $|\mathcal{O}_i| = m_i >1$. Define the auxiliary kernel
    $$\overline{M}_i(x,y) = 
    \begin{cases}
    \frac{1}{m_i}\alpha(x,y), &y\in\mathcal{O}_i,\ y \neq x,\\
    1 - \sum_{y\neq x}\overline{M}_i(x,y), & y = x,\\
    0, & \text{otherwise},
    \end{cases}$$
    which is stationary with respect to the distribution 
    $$\pi^{(i)} = \frac{1}{\pi(\mathcal{O}_i)}\bigg(\pi(x_1^{(i)}), \dots, \pi(x_{m_i}^{(i)})\bigg)$$
    on $\mathcal{O}_i$.

    By \cite{Liu1996}, the eigenvalues of $\overline{M}$ are given by 
    \begin{equation} \label{eq:7}
       \overline{\lambda_j} = 1 - \frac{j-2}{m_i} - \frac{1}{m_i\cdot \pi^{(i)}(x_{j-1}^{(i)})}\sum_{l = j-1}^{m_i} \pi^{(i)}(x^{(i)}_l). 
    \end{equation}
    Using \eqref{weyl} and the fact that $M_i = \frac{m_i}{m_i-1}\overline{M}_i - \frac{1}{m_i-1}I_{m_i}$,
    $$\lambda_2(M_i) = 1 - \frac{\pi(\mathcal{O}_i)}{(m_i-1) \pi(x^{(i)}_1)} \quad \text{and} \quad \lambda_{m_i}(M_i) = -\frac{\pi(x_{m_i}^{(i)})}{\pi(x_{m_i-1}^{(i)})(m_i-1)}.$$
    Then $\rho(M|_{S^\perp})$ is the largest absolute value among all such $\lambda_2(M_i)$ and $\lambda_{m_i}(M_i)$, since all eigenvalues that are equal to 1 has eigenvectors lying on $S$.
\end{proof}

\begin{remark} \label{G=I}
    In the case where every orbit $\mathcal{O}_i$ is of size 1, we would have $n$ orbits, and $G = M = I_n$. Trivially, $GPG = MPM$ in this case.
\end{remark}

If for every orbit, $|\mathcal{O}_i| > 2$ then a crude upper bound on $\theta$ would be 
$(m-2)/(m-1)$ where $m$ is the size of the largest orbit.

To prove this upper bound, for the first term, since $\pi(\mathcal{O}_i) \geq \pi(x_1^{(i)})$,
$$1 - \frac{\pi(\mathcal{O}_i)}{(m_i - 1)\pi(x_1^{(i)})} \leq 1 - \frac{1}{m_i-1} $$
and 
$$\frac{\pi(\mathcal{O}_i)}{(m_i - 1)\pi(x_1^{(i)})} - 1 \geq \frac{1}{m_i-1}-1.$$
For the second term, $\pi(x_{m_i-1}^{(i)}) \geq \pi(x_{m_i}^{(i)}),$
$$\frac{\pi(x_{m_i}^{(i)})}{\pi(x_{m_i-1}^{(i)})(m_i-1)} \leq \frac{1}{m_i-1}.$$
The first term's upper bound dominates and is thus an upper bound of $\theta$.

For $\varepsilon > 0$, we define a time $t(\varepsilon)$ as follows:
\begin{align*}
    t(\varepsilon) := \inf\{k \in \mathbb{N};~ \max_i|\lambda_i(M^k P M^k) - \lambda_i(GPG)| \leq \varepsilon\}.
\end{align*}
This $t(\varepsilon)$ can be intuitively understood as the time it takes to approximate $G$ using $M^{t(\varepsilon)}$. Using Proposition \ref{SLEMM^kPM^k}, we thus have
\begin{align*}
    t(\varepsilon) \leq \max\bigg\{\dfrac{\ln(4\rho(P)/\varepsilon)}{\ln(1/\theta)},\dfrac{\ln(2\rho(P)/\varepsilon)}{2\ln(1/\theta)}\bigg\} =: \overline{t}(\varepsilon,\theta,P) = \overline{t}.
\end{align*}
As a result, a heuristic is that we use $M^{\overline{t}}PM^{\overline{t}}$ to approximate $GPG$.

In the last part of this subsection, we make a remark that there is no strict ordering between $\rho(BPB)$ and $\rho(MPM)$ in general. While it is well-known by the results of \cite{Peskun1973} that the Metropolis-Hastings proposal is the best in its family under certain metrics, such optimality results is no longer guaranteed when we consider their multiplicative sandwich $BPB$ and $MPM$. 

\subsection{The Metropolis-Hastings orbit sampler with one orbit}
We aim to provide a concrete example where the Metropolis-Hastings orbit sandwich $M^k P M^k$ outperforms the kernel $P$ by an exponential order. 

Let $\psi$ be the uniform distribution on $\mathcal{X} = \llbracket n \rrbracket$. Consider the lazy random walk sampler
$$P =
    \begin{pmatrix}
    \frac{1}{2} & \frac{1}{2} & 0 & 0 & \cdots & 0 & 0 \\
    \frac{1}{4} & \frac{1}{2} & \frac{1}{4} & 0 & \cdots & 0 & 0 \\
    0 & \frac{1}{4} & \frac{1}{2} & \frac{1}{4} & \cdots & 0 & 0 \\
    \vdots & & \ddots & \ddots & \ddots & & \vdots \\
    0 & 0 & \cdots & \frac{1}{4} & \frac{1}{2} & \frac{1}{4} & 0 \\
    0 & 0 & \cdots & 0 & \frac{1}{4} & \frac{1}{2} & \frac{1}{4} \\
    0 & 0 & \cdots & 0 & 0 & \frac{1}{2} & \frac{1}{2}
    \end{pmatrix}.
$$
It can be verified that it is a reversible sampler with respect to the uniform distribution. Furthermore, its eigenvalues are 
$$\lambda_m = \frac{1}{2} + \frac{1}{2}\cos \bigg( \frac{(m-1)\pi}{n-1}\bigg)$$
with the associated eigenvectors
$$
v_m(i) = \begin{cases}
    1, & m = 1,\\
    \cos\big(\frac{(m+1)\pi}{n-1}(i-1)\big), & m = 2, \dots, n-1,\\
    (-1)^{i-1}, & m = n.
\end{cases}
$$
Note that $\pi$ in the results refer to the constant, and not the stationary distribution which we denote by $\psi$ to avoid confusion. We shall still use $\Pi$ to denote the matrix with rows all equal to $\psi$.

Hence, $\rho(P)$ is of the order $1-\Theta(n^{-2})$ which indicates a diffusive-type mixing. 

Now suppose that our choice of group $\mathcal{G}$ admits only a single group. The Gibbs kernel $G = \Pi$ and the MH kernel $M$ is of the form 
$$
M(x,y) = \begin{cases}
    \frac{1}{n-1}, & \text{if } x \neq y,\\
    0, & \text{otherwise}.
\end{cases}
$$
Then $\theta$, the largest absolute eigenvalue in $\mathrm{spec}(M|_{S^\perp})$ is $(n-1)^{-1}$. 

By Proposition \ref{SLEMM^kPM^k}, the convergence rate of $M^kPM^k$ to $GPG = \Pi$ is given by $\rho(P)(2\theta^k + \theta^{2k})$ which is of the order $\Theta(n^{-k})\rho(P)$. Compared to both $\rho(M)$ and $\rho(P)$, we see exponential improvement by using the sandwich $M^kPM^k$. In fact, merely using $k = 1$ (i.e. $MPM$) leads to a constant independent of $n$ order in relaxation time since $\rho(P) (2\theta^k + \theta^{2k}) = \Theta(n^{-1}) = O(1)$.\\

\begin{remark}
    Even in the case where $\pi$ is non-uniform, one can use the results in Proposition \ref{thetadef} to find an upper bound on $\theta$. In the case of a single orbit, 
    $$\theta = \max\bigg\{ \bigg|1-\frac{1}{(n-1)\pi(x_1)}\bigg|, \frac{\pi(x_n)}{(n-1)\pi(x_{n-1})}\bigg\},$$
    where $\pi(x_1) \geq \dots \geq \pi(x_n)$. If $\pi(x)$ is in the order of $\Theta(n^{-1})$ for all $x \in \mathcal{X}$, one can still expect exponential improvement in SLEM when using $M^kPM^k$.
\end{remark}

\section{Comparison of asymptotic variance} \label{section:var}

In this section, we investigate and compare the asymptotic variance of the samplers $P$, $BPB$, $MPM$ and $GPG$. 

Let $P\in \mathcal{S}(\pi)$ be an ergodic transition matrix. We write as per \cite{Bremaud1999}, Ch. 6,
$$Z(P) = (I - (P - \Pi))^{-1}$$ 
to be the fundamental matrix of $P$, where $\Pi$ is the matrix with each row as $\pi$. 

Then, the asymptotic variance of $f\in \ell^2_0(\pi)$, for any initial distribution $\mu$ is given by \cite{Bremaud1999}, Theorem 6.5, as
\begin{equation} \label{asympvar}
    v(f,P) := \lim_{n\rightarrow\infty}\frac{1}{n}\text{Var}\bigg( \sum_{i=1}^n f(X_i) \bigg) = 2\langle f,Z(P)f\rangle_\pi - \langle f,f\rangle_\pi.
\end{equation}

An equivalent characterisation given by \cite{sherlock2025reversiblemarkovchainsvariational} for $v(f, P)$ when $P \in \mathcal{L}(\pi)$ is 
\begin{equation} \label{asymptvar2}
    v(f,P) = \sup_{h\in \ell^2_0(\pi)} 4\langle f,h\rangle_\pi - 2\langle (I-P)h,h \rangle_\pi - \langle f,f\rangle_\pi.
\end{equation}

The worst-case asymptotic variance, given by \cite{Frigessi1993} is
\begin{equation} \label{worsecasevar}
    V(P) := \sup_{f \in \ell^2_0(\pi),\ \|f\|_\pi = 1} v(f,P) = \frac{1+\lambda_2(P)}{1-\lambda_2(P)}.
\end{equation}

\begin{prop} \label{varGPG/P}
Let $P\in\mathcal{L}(\pi)$ be an ergodic transition matrix on $\mathcal X$, and let $G$ be the Gibbs kernel induced by any group action $\mathcal G$ on $\mathcal{X}$. Write
\[
S=\operatorname{Ran}(G),\qquad f_\perp=(I-G)f.
\]
Then the following statements hold.
\begin{enumerate}[(i)]
    \item For every $f\in\ell^2_0(\pi)$,
    \begin{equation}
        \label{eq:additive-variance-bound}
        v(f,GPG)
        \leq
        v(f,P)+2\|(I-G)f\|_\pi^2.
    \end{equation}

    \item In particular, for every $G$-invariant observable
    $f\in S\cap\ell^2_0(\pi)$, that is, $f$ satisfies $Gf = f$, 
    \[
        v(f,GPG)\leq v(f,P).
    \]
    Equality holds if and only if $Z(P)f\in S$. Equivalently,
    \[
        \left((I-P)|_{\ell^2_0(\pi)}\right)^{-1}f\in S.
    \]
\end{enumerate}
\end{prop}

\begin{proof}
Every $h\in\ell^2_0(\pi)$ admits the orthogonal decomposition
\[h=u+w, \qquad u=Gh\in S\cap\ell^2_0(\pi), \qquad w=(I-G)h\in S^\perp\cap\ell^2_0(\pi).\]

Since $Gh=u$, we have $GPGh=GPu$. Therefore, using the self-adjointness of $G$,
\begin{align}
    \langle(I-GPG)h,h\rangle_\pi &= \|u+w\|_\pi^2-\langle GPu,u+w\rangle_\pi \notag\\
    &= \|u\|_\pi^2+\|w\|_\pi^2 -\langle Pu,G(u+w)\rangle_\pi \notag\\
    &= \|u\|_\pi^2-\langle Pu,u\rangle_\pi+\|w\|_\pi^2 \notag\\
    &= \langle(I-P)u,u\rangle_\pi+\|w\|_\pi^2.
    \label{eq:GPG-dirichlet-decomposition}
\end{align}

Also, writing $f=Gf+f_\perp$, orthogonality gives
\begin{align*}
\langle f,h\rangle_\pi &= \langle f,u\rangle_\pi+\langle f_\perp,w\rangle_\pi.
\end{align*}

Using the variational representation \eqref{asymptvar2} and
\eqref{eq:GPG-dirichlet-decomposition}, we obtain
\begin{align*}
v(f,GPG)
&=
\sup_{h\in\ell^2_0(\pi)}
\left\{
4\langle f,h\rangle_\pi
-2\langle(I-GPG)h,h\rangle_\pi
\right\}
-\|f\|_\pi^2\\
&=
\sup_{\substack{u\in S\cap\ell^2_0(\pi)\\
                 w\in S^\perp\cap\ell^2_0(\pi)}}
\bigg\{
4\left(
\langle f,u\rangle_\pi
+\langle f_\perp,w\rangle_\pi
\right)
-2\left(
\langle(I-P)u,u\rangle_\pi+\|w\|_\pi^2
\right)
\bigg\}
-\|f\|_\pi^2\\
&=
\sup_{\substack{u\in S\cap\ell^2_0(\pi)\\
                 w\in S^\perp\cap\ell^2_0(\pi)}}
\bigg\{
4\langle f,u\rangle_\pi
+4\langle f_\perp,w\rangle_\pi
-2\langle(I-P)u,u\rangle_\pi
-2\|w\|_\pi^2
\bigg\}
-\|f\|_\pi^2.
\end{align*}
The objective is a sum of a function of $u$ and a function of $w$,
so the supremum separates. Moreover,
\begin{align*}
4\langle f_\perp,w\rangle_\pi-2\|w\|_\pi^2
&=
2\|f_\perp\|_\pi^2
-2\|w-f_\perp\|_\pi^2.
\end{align*}
Since $f_\perp\in S^\perp\cap\ell^2_0(\pi)$, it follows that
\[
\sup_{w\in S^\perp\cap\ell^2_0(\pi)}
\left\{
4\langle f_\perp,w\rangle_\pi-2\|w\|_\pi^2
\right\}
=
2\|f_\perp\|_\pi^2,
\]
with the supremum attained at $w=f_\perp$. Hence
\begin{align*}
v(f,GPG)
={}&
\sup_{u\in S\cap\ell^2_0(\pi)}
\bigg\{
4\langle f,u\rangle_\pi
-2\langle(I-P)u,u\rangle_\pi
\bigg\}
+2\|f_\perp\|_\pi^2-\|f\|_\pi^2\\
\leq{}&
\sup_{h\in\ell^2_0(\pi)}
\bigg\{
4\langle f,h\rangle_\pi
-2\langle(I-P)h,h\rangle_\pi
\bigg\}
+2\|f_\perp\|_\pi^2-\|f\|_\pi^2\\
={}&
v(f,P)+2\|f_\perp\|_\pi^2.
\end{align*}
This proves part~(i). Part~(ii) follows immediately upon taking
$f_\perp=0$.

It remains to establish the equality condition in part~(ii). Since
$P$ is ergodic, the eigenvalue $1$ of $P$ is simple. Therefore, $1$ is not an eigenvalue
of the restriction of $P$ to $\ell^2_0(\pi)$, and hence
\[
(I-P)|_{\ell^2_0(\pi)}
\]
is invertible.

For fixed $f\in\ell^2_0(\pi)$, the variational problem defining $v(f,P)$
has the unique maximiser
\[
h_*
=
\left((I-P)|_{\ell^2_0(\pi)}\right)^{-1}f.
\]
Indeed, the first-order optimality condition is
\[
(I-P)h_*=f,
\]
and uniqueness follows because $I-P$ is positive definite on
$\ell^2_0(\pi)$.

The variational problem defining $v(f,GPG)$ for $f\in S$ is precisely
the same problem restricted to $S\cap\ell^2_0(\pi)$. Consequently,
\[
v(f,GPG)=v(f,P)
\]
if and only if the unique unrestricted maximiser $h_*$ belongs to
$S$.

Finally, since $\Pi f=0$ for every $f\in\ell^2_0(\pi)$,
\[
(I-P+\Pi)|_{\ell^2_0(\pi)}=(I-P)|_{\ell^2_0(\pi)}.
\]
It follows that
\[
h_*
=
\left((I-P)|_{\ell^2_0(\pi)}\right)^{-1}f
=
Z(P)f.
\]
Thus equality holds if and only if $Z(P)f\in S$.
\end{proof}

In general, such ordering of asymptotic variance between $MPM$ and $BPB$ against $P$ does not exist. However, we can turn to the worst-case asymptotic variance, and comparisons can be made between all of them assuming certain conditions. 

\begin{prop} \label{ordervar}
    For any $\pi$-stationary and reversible ergodic sampler $P$, the inequalities 
    \begin{align*}
        V(GPG) &\leq V(MPM) \leq V(P), \\
        V(GPG) &\leq V(BPB) \leq V(P).
    \end{align*}
    hold if $P$ admits non-negative spectra. 
\end{prop}

\begin{proof}
    If $\mathrm{spec}(P)$ lie between $[0,1]$, then $\rho(P) = \lambda_2(P)$, and similarly for all the multiplicative group-averaged chains $GPG$, $MPM$ and $BPB$. From Proposition \ref{SLEMGPG}, Proposition \ref{SLEMMPM} and Proposition \ref{SLEMdiffGPG/MPM}, 
    \begin{align*}
        \lambda_2(GPG) &\leq \lambda_2(MPM) \leq \lambda_2(P),\\
        \lambda_2(GPG) &\leq \lambda_2(BPB) \leq \lambda_2(P).
    \end{align*} 
    Since $V(P)$ increases with $\lambda_2(P)$ when $\lambda_2 \in [0,1]$, the inequality follows.
\end{proof}

\section{Pythagorean identity and comparison of one-step KL divergence to stationarity} \label{section:KL}
Let $\pi \in \mathcal{P}(\mathcal{X})$ be a probability mass. For $P, Q \in \mathcal{L}$, the KL divergence of $P$ from $Q$ is defined as 
\begin{equation} \label{defKL}
    D_{KL}^\pi (P \| Q) := \sum_{x,y \in \mathcal{X}} \pi(x) P(x,y) \log\bigg(\frac{P(x,y)}{Q(x,y)}\bigg),
\end{equation}
where by convention we take $0 \log (0/a) := 0$ for $a \in [0,1]$.
With a chosen group $\mathcal{G}$ and its corresponding Gibbs orbit kernel $G$, let $\mathbf{G}$ (resp.~$\mathbf{M}, \mathbf{B}$) be the set of invariant samplers under $GPG$ (resp.~$MPM, BPB$). Formally,
\begin{align*}
    \mathbf{G} &= \mathbf{G}(\mathcal{G},\pi) := \{P \in \mathcal{S}(\pi): GPG = P\},\\
    \mathbf{M} &= \mathbf{M}(\mathcal{G},\pi) := \{P \in \mathcal{S}(\pi): MPM = P\},\\
    \mathbf{B} &= \mathbf{B}(\mathcal{G},\pi) := \{P \in \mathcal{S}(\pi): BPB = P\}.
\end{align*}

Under most circumstances, the invariant sets $\mathbf{G} = \mathbf{M} = \mathbf{B}$ coincide. \\

\begin{prop}\label{prop:mathbfGMB}
    For a fixed group action $\mathcal{G}$ and the invariant sets $\mathbf{G}, \mathbf{M}, \mathbf{B}$, if each block in $G, M$ and $B$ is aperiodic then 
    $$\mathbf{G} = \mathbf{M} = \mathbf{B}.$$
\end{prop}

\begin{proof}
    Suppose if $P \in \mathbf{G}$, then 
    $$MPM = MGPGM = GPG = P,$$
    noting that $GM = MG = G$ by Lemma \ref{GM=MG}. 

    If instead $P \in \mathbf{M}$, and given that $M$ is aperiodic, Proposition \ref{limMB} gives $M^t \rightarrow G$ as $t \rightarrow \infty$ and so
    $$P = MPM = M^2 P M^2 = \dots = GPG.$$
    Hence, $P \in \mathbf{G}$ if it is in $\mathbf{M}$. 
    
    By replacing $M$ with $B$, we then also have that $\mathbf{B} = \mathbf{M} = \mathbf{G}$.
\end{proof}

Next, we give a characterisation of the $\pi$-stationary kernels on $\mathcal{X}$ that lie in $\mathbf{G}$ for a given $\mathcal{G}.$\\

\begin{prop}
    Let $\mathcal{G}$ define orbits $(\mathcal{O}_i)_{i=1}^k$. Then $Q \in \mathcal{S}(\pi)$ lies in $\mathbf{G}$ if and only if, for $x \in \mathcal{O}_i$ and $y \in \mathcal{O}_j$, 
    $$Q(x,y) = c_{ij} \frac{\pi(y)}{\pi(\mathcal{O}_j)},$$
    for some coefficients $c_{ij} \geq 0$ satisfying 
    \begin{equation} \label{c_ij}
        \sum_{i=1}^k c_{ij} = 1 \quad \text{and} \quad \sum_{i=1}^k \pi(\mathcal{O}_i)c_{ij} = \pi(\mathcal{O}_j).
    \end{equation}
\end{prop}

\begin{proof}
    Partition $Q$ into orbit blocks $Q_{ij} \in \mathbb{R}^{|\mathcal{O}_i| \times |\mathcal{O}_j|}$, and write $G = \text{diag}(G_1, \dots, G_k)$ with off-diagonal blocks as zeros. Define
    $$\mu_i(x) := \frac{\pi(x)}{\pi(\mathcal{O}_i)},$$
    and let $\mathbf{1}_i$ be the $|\mathcal{O}_i| \times 1$ column vector of 1's. 
    
    Then each $G_i$ is the Gibbs orbit kernel on $\mathcal{O}_i$, and $G_i = \mathbf{1}_i \mu_i^T. $

    First, suppose that $Q \in \mathbf{G}$, that is $GQG = Q.$ Blockwise, 
    $$Q_{ij} = G_i Q_{ij} G_j = (\mathbf{1}_i \mu_i^T) Q_{ij} (\mathbf{1}_j \mu_j^T).$$
    Set $c_{ij} = \mu_i^T Q_{ij} \mathbf{1}_j$, and by the fact that $Q$ is $\pi$-stationary, we have that $c_{ij}$ satisfies \eqref{c_ij}.

    Now suppose $Q_{ij} = \mathbf{1}_i c_{ij} \mu_j^T$, with $c_{ij}$ satisfying \eqref{c_ij}. Then for any $(i,j)$ pair, the block
    $$(GQG)_{ij} = G_i Q_{ij} G_j = G_i (\mathbf{1} c_{ij} \mu_j^T) G_j.$$
    Since each $G_i$ is a kernel on $\mathcal{O}_i$ with stationary distribution $\mu_i$,
    $$(GQG)_{ij} = \mathbf{1} c_{ij} \mu_j^T = Q_{ij}.$$
\end{proof}

With the results above, for $Q$ to be in $\mathbf{G}$, all its rows within an orbit $\mathcal{O}_i$ must be identical. Furthermore, the columns in each $(i,j)$ block must be proportional to the stationary weights $\pi(y)$. 

We now show that for any $P \in \mathcal{L}(\pi)$, we have the following Pythagorean identity.\\ 

\begin{prop}\label{prop:PythidGPG} 
    Let $G$ be a Gibbs orbit kernel, with the orbits $(\mathcal{O}_i)_{i=1}^k$. For $P \in \mathcal{S}(\pi)$ and $Q \in \mathbf{G}$, 
    \begin{equation} \label{KLineqGPG}
        D_{KL}^\pi(P \| Q) = D_{KL}^\pi(P \| GPG) + D_{KL}^\pi(GPG \| Q).
    \end{equation}
    In particular, this implies that $GPG$ is the unique projection of $P$ onto $\mathbf{G}$ under the KL divergence, that is,
    \begin{align*}
        D_{KL}^\pi(P \| GPG) = \min_{Q \in \mathbf{G}} D_{KL}^\pi(P \| Q).
    \end{align*}
    By replacing $P$ above with either $MPM$ or $BPB$ and noting Lemma \ref{GM=MG}, we see that
    \begin{align*}
        D_{KL}^\pi(MPM \| Q) &= D_{KL}^\pi(MPM \| GPG) + D_{KL}^\pi(GPG \| Q), \\
        D_{KL}^\pi(BPB \| Q) &= D_{KL}^\pi(BPB \| GPG) + D_{KL}^\pi(GPG \| Q).
    \end{align*}
    In other words, $GPG$ is also the unique projection of either $MPM$ or $BPB$ onto $\mathbf{G}$ under the KL divergence.

    By specializing into $Q = \Pi$, we have
    \begin{align*}
        D_{KL}^\pi(P \| \Pi) &\geq D_{KL}^\pi(GPG \| \Pi), ~D_{KL}^\pi(MPM \| \Pi) \geq D_{KL}^\pi(GPG \| \Pi), \\
        D_{KL}^\pi(BPB \| \Pi) &\geq D_{KL}^\pi(GPG \| \Pi).
    \end{align*}
\end{prop}

\begin{proof}
    Notice that
    \begin{align*}
         D_{KL}^\pi(P \| Q) &= D_{KL}^\pi(P \| GPG) + D_{KL}^\pi(GPG \| Q)\\ 
         &+ \sum_{x,y \in \mathcal{X}} \pi(x) \big( P(x,y) - GPG(x,y)\big) \log \bigg(\frac{GPG(x,y)}{Q(x,y)}\bigg),
    \end{align*}
    and so it suffice to show that the last term on the right is 0.

    Since $Q \in \mathbf{G}$, using \eqref{GPGmatrix}, the sum can be written as 
    \begin{align*}
        &\sum_{x,y \in \mathcal{X}} \pi(x) \big( P(x,y) - GPG(x,y)\big) \log \bigg(\frac{GPG(x,y)}{Q(x,y)}\bigg)\\
        = &\sum_{i,j = 1}^k \sum_{\substack{x \in \mathcal{O}_i\\ y \in \mathcal{O}_j}} \pi(x)\big( P(x,y) - GPG(x,y)\big) \log \Bigg(\frac{\sum_{\substack{z \in \mathcal{O}_i\\ w \in \mathcal{O}_j}} \pi(z)P(z,w)}{\sum_{\substack{z \in \mathcal{O}_i\\ w \in \mathcal{O}_j}} \pi(z)Q(z,w)} \Bigg)\\
        &= \sum_{i,j = 1}^k \log \Bigg(\frac{\sum_{\substack{z \in \mathcal{O}_i\\ w \in \mathcal{O}_j}} \pi(z)P(z,w)}{\sum_{\substack{z \in \mathcal{O}_i\\ w \in \mathcal{O}_j}} \pi(z)Q(z,w)} \Bigg) \sum_{\substack{x \in \mathcal{O}_i\\ y \in \mathcal{O}_j}} \pi(x)\big( P(x,y) - GPG(x,y)\big).
    \end{align*}

    By Proposition \ref{projchain_equal}, the inner sum must be 0 for any $i,j \in \llbracket k \rrbracket.$
\end{proof}

\begin{figure}[ht]
  \centering
  \begin{tikzpicture}[>=Latex, every node/.style={font=\small},
  arrowmid/.style={postaction={decorate,decoration={
    markings,
    mark=at position 0.55 with {\arrow[scale=1.8]{Latex[length=4pt,width=3pt]}}
  }}}
]

  %--- Parameters ---
  \def\ax{4.6}
  \def\ay{2.7}

  %--- Draw the common space (ellipse) ---
  \def\ay{1.2} % instead of 2.7
\draw[thick] (0,0) ellipse ({\ax} and {\ay});
  \node[below=2pt] at (0,-\ay) {$\mathbf{B}=\mathbf{M}=\mathbf{G}$};

  %--- Points inside ---
  \coordinate (GPG) at (0,-0.35);
  \fill (GPG) circle (1.3pt);
  \node[below right=2pt] at (GPG) {$GPG$};

  %--- Down arrow to Π ---
  \coordinate (Pi) at (0,-1.0);
  \fill (Pi) circle (1pt);
  \draw[arrowmid] (Pi) -- (GPG);
  \node[left=0.5pt] at (Pi) {$\Pi$};

  %--- Exterior points ---
  \coordinate (P)   at (0,5.1);
  \coordinate (MPM) at (-4.2,3.2);
  \coordinate (BPB) at ( 4.2,3.2);

  \fill (P) circle (1.2pt) node[above=3pt] {$P$};
  \fill (MPM) circle (1.2pt) node[left=3pt] {$MPM$};
  \fill (BPB) circle (1.2pt) node[right=3pt] {$BPB$};

  %--- Control point for P (curved path) ---
  \coordinate (CPp) at (0,3.8);

  %--- Projection curves ---
  \draw[arrowmid] (GPG) .. controls (CPp) .. (P);
  \draw[arrowmid] (GPG) -- (MPM);
  \draw[arrowmid] (GPG) -- (BPB);

  %--- Horizontal line at GPG ---
  \draw[thick] ($(GPG)+(-0.8,0)$) -- ($(GPG)+(0.8,0)$);

  %--- Perpendicular from P to GPG ---
  \draw[dashed] (P) -- (GPG);
  \draw ($(GPG)+(-0.2,0)$) -- ++(0,0.2) -- ++(0.2,0); % right-angle symbol at GPG

  %===============================
  % Transparent extensions + perpendicular markers
  %===============================

  % Extend lines slightly past GPG (into ellipse)
  \coordinate (ExtMPM) at ($ (GPG)!-0.25!(MPM) $);
  \coordinate (ExtBPB) at ($ (GPG)!-0.25!(BPB) $);
  \draw[gray!25,dashed] (GPG) -- (ExtMPM);
  \draw[gray!25,dashed] (GPG) -- (ExtBPB);

  % Compute directions
  \pgfmathanglebetweenpoints{\pgfpointanchor{GPG}{center}}{\pgfpointanchor{MPM}{center}}
  \let\angGtoMPM\pgfmathresult
  \pgfmathanglebetweenpoints{\pgfpointanchor{GPG}{center}}{\pgfpointanchor{BPB}{center}}
  \let\angGtoBPB\pgfmathresult

  % Positions for perpendicular markers (just before GPG)
  \coordinate (MPMperp) at ($ (GPG)!0.045!(MPM) $);
  \coordinate (BPBperp) at ($ (GPG)!0.035!(BPB) $);

  % Size and offset of perpendicular marks
  \def\perpsizeM{0.20}
  \def\offsetM{0.20}

  \def\perpsizeB{0.15}
  \def\offsetB{0.15}

  % MPM side (inward)
  \begin{scope}[shift={(MPMperp)}, rotate=\angGtoMPM]
    \draw[line width=0.4pt, shift={(0,\offsetM)}] (0,0) -- (\perpsizeM,0) -- (\perpsizeM,-\perpsizeM);
  \end{scope}

    % BPB side (flipped outward, upright, rotated 180°)
    \begin{scope}[shift={(BPBperp)}, rotate=\angGtoBPB+180]
      \draw[line width=0.4pt, shift={(0,\offsetB)}] (0,0) -- (-\perpsizeB,0) -- (-\perpsizeB,-\perpsizeB);
    \end{scope}

\end{tikzpicture}
  \caption{Visualisation of $\mathbf{G}$ and the projections of various samplers onto $\mathbf{G}$ under the assumptions of Proposition \ref{prop:mathbfGMB} and \ref{prop:PythidGPG}.}
  \label{fig:projection-diagram}
\end{figure}

While Proposition \ref{prop:PythidGPG} establishes the Pythagorean identity for the Gibbs orbit kernel, the same relationship does not generally hold when $G$ is replaced with $M$ or $B$. In fact, for any $Q \in \mathbf{G}$, there is no uniform ordering between $ D_{KL}^\pi(P \| Q)$ and $ D_{KL}^\pi(P \| MPM) + D_{KL}^\pi(MPM \| Q)$ (resp. $BPB$).

The following counterexample illustrates this.

Let $\pi = (0.1, 0.2, 0.3, 0.4)$, and choose $\mathcal{G}$ such that the orbits are $(1,2)$ and $(3,4)$. Define $G$ to be the associated Gibbs orbit kernel, and let $P$ denote the usual Metropolis–Hastings kernel for $\pi$:
$$P = \begin{pmatrix}
    0 & 1/3 & 1/3 & 1/3\\
    1/6 & 1/6 & 1/3 & 1/3\\
    1/9 & 2/9 & 1/3 & 1/3\\
    1/12 & 1/6 & 1/4 & 1/2
    \end{pmatrix},$$
and set $Q = GPG \in \mathbf{G}$. 

Direct computation yields
$$D_{KL}^\pi(P \| Q) = 0.0301 < 0.03702 = D_{KL}^\pi(P \| MPM) + D_{KL}^\pi(MPM \| Q),$$
showing that the Pythagorean decomposition fails.

However, if we instead consider the lazified kernel $P_0 = \tfrac{1}{2}(I + P)$, the inequality reverses:
$$D_{KL}^\pi(P_0 \| Q) = 0.29026 > 0.21660 = D_{KL}^\pi(P_0 \| MP_0M) + D_{KL}^\pi(MP_0M \| Q).$$
Thus, even the direction of the inequality depends on the particular form of the transition kernel.
This highlights that the exact orthogonality property is unique to the Gibbs sampler $GPG$.

Now, even though the Pythagorean identity does not generally hold for $M$ or $B$, these kernels still act as KL-contractive steps towards $\mathbf{G}$.\\

\begin{prop} \label{BP/PB}
    Let $M$ and $B$ be the Barker and MH orbit sampler with the orbits $(\mathcal{O}_i)_{i=1}^k$ respectively. For $P \in \mathcal{S}(\pi)$ and $Q \in \mathbf{G}$, we have the inequalities
    \begin{align*}
        D_{KL}^\pi(P \| Q) \geq D_{KL}^\pi(PM \| Q),\  D_{KL}^\pi(P \| Q) \geq D_{KL}^\pi(PB \| Q),\\
        D_{KL}^\pi(P \| Q) \geq D_{KL}^\pi(MP \| Q),\  D_{KL}^\pi(P \| Q) \geq D_{KL}^\pi(BP \| Q).
    \end{align*}
    These inequalities can be interpreted as an analogue of the data-processing inequality in our context.
\end{prop}

\begin{proof}
    Consider 
    $$D_{KL}^\pi(P \| Q) - D_{KL}^\pi(PB \| Q) = \sum_{i,j = 1}^k \sum_{\substack{x \in \mathcal{O}_i \\ y \in \mathcal{O}_j}} \pi(x) \bigg( P(x,y)\log\bigg(\frac{P(x,y)}{Q(x,y)}\bigg) - PB(x,y)\log\bigg(\frac{PB(x,y)}{Q(x,y)}\bigg) \bigg).$$

    For any  fixed $x \in \mathcal{O}_i$ and $y \in \mathcal{O}_j$, the log-sum inequality gives
    \begin{align*}
        \sum_{z \in \mathcal{X}} P(x,z)B(z,y) \log\bigg(\frac{P(x,z)B(z,y)}{Q(x,z) B(z,y)}\bigg) \geq PB(x,y)\log \bigg(\frac{PB(x,y)}{Q(x,y)}\bigg).
    \end{align*}
    Summing across all possible $y \in \mathcal{X}$ on both sides,
    $$\sum_{z \in \mathcal{X}} P(x,z) \log\bigg(\frac{P(x,z)}{Q(x,z)}\bigg) \geq \sum_{y \in \mathcal{X}} PB(x,y)\log \bigg(\frac{PB(x,y)}{Q(x,y)}\bigg).$$

    Then, by multiplying $\pi(x)$ and summing up over all possible $x$, 
    \begin{align*}
        D_{KL}^\pi (P \| Q) &= \sum_{x,y \in \mathcal{X}} \pi(x) P(x,y) \log\bigg(\frac{P(x,y)}{Q(x,y)}\bigg)\\
        &\geq \sum_{x, y \in \mathcal{X}} \pi(x) PB(x,y)\log \bigg(\frac{PB(x,y)}{Q(x,y)}\bigg) = D_{KL}^\pi (PB \| Q)
    \end{align*}

    With the bisection property $D^\pi_{KL}(P \| Q) = D^\pi_{KL}(P^* \| Q^*)$ shown in \cite{Choi_2024} Theorem 3.1, the other inequality follows from
    \begin{align*}
        D_{KL}^\pi (P \| Q) &= D_{KL}^\pi (P^* \| Q^*)\\
        &\geq D_{KL}^\pi (P^* B \| Q^*)\\ 
        &= D_{KL}^\pi (BP \| Q).
    \end{align*}

    By replacing $B$ with $M$ above, one can obtain the other two inequalities.
\end{proof}

By collecting the previous two results, we arrive at the following Corollary:

\begin{corollary}
    For $P \in \mathcal{S}(\pi)$ and $Q \in \mathbf{G}$, 
    \begin{align*}
        &D_{KL}^\pi(P \| Q) \geq D_{KL}^\pi(MPM \| Q) \geq D_{KL}^\pi(M^2 P M^2 \| Q) \geq \ldots \geq D_{KL}^\pi(GPG \| Q),\\
        &D_{KL}^\pi(P \| Q) \geq D_{KL}^\pi(BPB \| Q) \geq D_{KL}^\pi(B^2 P B^2 \| Q) \geq \ldots \geq D_{KL}^\pi(GPG \| Q),
    \end{align*}
\end{corollary}

\begin{proof}
    From Proposition \ref{BP/PB}, it follows that 
    $$D_{KL}^\pi(P \| Q) \geq D_{KL}^\pi(PM \| Q) \geq D_{KL}^\pi(MPM \| Q),$$
    and so inductively, for all $n \geq 1$,
    $$D_{KL}^\pi(M^n P M^n\| Q) \geq D_{KL}^\pi(M^{n+1} P M^{n+1}\| Q).$$

    Furthermore, for any $n \geq 1$, Proposition \ref{prop:PythidGPG} gives 
    $$D_{KL}^\pi(M^n P M^n \| Q) \geq D_{KL}^\pi(G M^n P M^n G\| Q) = D_{KL}^\pi(GPG \| Q),$$
    since $MG = GM = G$ by Proposition \ref{GM=MG}.

    The proof is identical for the case of $BPB$.
\end{proof}

\section{Orbit-space reduction and construction of samplers for a fixed group action} \label{section:bestP}
This section investigates the effect of fixing the orbit partition induced by a group action and varying the base sampler $P$. We first show that the spectral and Kullback--Leibler properties of $GPG$ reduce exactly to those of the projection chain $\overline P$ on the orbit space. We then use this reduction to construct a particular family of samplers with favourable behaviour when the stationary mass is concentrated on one orbit.

Let $\widetilde{P}$ be a sampler on the orbit space $(\mathcal{O}_i)^k_{i=1}$, that is stationary and reversible with respect to the distribution $\overline{\pi} = (\pi(\mathcal{O}_1), \dots, \pi(\mathcal{O}_k)).$ Define the orbit-average sampler $Q_{\widetilde{P}} = Q_{\widetilde{P}}(\mathcal{G})$ on $\llbracket n \rrbracket$ as
\begin{equation} \label{liftingP}
    Q_{\widetilde{P}}(x,y) := \displaystyle \widetilde{P}(i,j)\,\frac{\pi(y)}{\pi(\mathcal O_j)},\quad \text{for } x\in\mathcal O_i,\ y\in\mathcal O_j.
\end{equation}
It can be verified that $Q_{\widetilde{P}}$ is both stationary and reversible with respect to $\pi$.\\

Furthermore, one can define an isometry $U: \mathbb{R}^k \to S$,
\begin{equation} \label{U}
    (Uf)(x) = f(i) \text{ for } x \in \mathcal{O}_i.
\end{equation}

Its adjoint $U^*$ satisfying $\langle Uf, g \rangle_\pi = \langle f, U^*g\rangle_{\overline{\pi}}$ for $f \in \ell^2(\overline{\pi}),\ g \in \ell^2(\pi)$ is given by 
\begin{equation} \label{U*}
    (U^*g)(i) = \frac{1}{\pi(\mathcal{O}_i)} \sum_{x \in \mathcal{O}_i} \pi(x) g(x).
\end{equation}

It then holds that $U^* U = I$ on $\ell^2(\overline{\pi})$ and $UU^* = G$. This isometry is the key connection between the two state space $\llbracket n \rrbracket$ and $(\mathcal{O}_i)^k_{i=1}$.\\

\begin{prop} \label{spectral orbit = sandwich}
    For any non-trivial group action $\mathcal{G}$ with $k < n$ orbits, the non-trivial spectrum of $G Q_{\widetilde{P}} G$ is exactly that of $\widetilde{P}$. That is,
    $$\normalfont \mathrm{spec}(G Q_{\widetilde{P}} G)  = \mathrm{spec}(\widetilde{P}) \cup \{0\}.$$
\end{prop}

\begin{proof}
     Let $U$ and $U^*$ be defined as per \eqref{U} and \eqref{U*} respectively. For any $g \in \ell^2(\overline{\pi})$, 
    \begin{align*}
        (U^* Q_{\widetilde{P}} Ug)(i) &= \frac{1}{\pi(\mathcal{O}_i)} \sum_{x \in \mathcal{O}_i} \pi(x) (Q_{\widetilde{P}} Ug)(x)\\
        &= \frac{1}{\pi(\mathcal{O}_i)} \sum_{x \in \mathcal{O}_i} \bigg[\pi(x) \sum_{j=1}^k \bigg(g(j) \sum_{y \in \mathcal{O}_j}Q_{\widetilde{P}}(x,y)\bigg) \bigg]\\
        &= \frac{1}{\pi(\mathcal{O}_i)} \sum_{x \in \mathcal{O}_i} \bigg( \pi(x) \sum_{j=1}^k g(j) \widetilde{P}(i,j)\bigg)\\
        &= \sum_{j=1}^k g(j) \widetilde{P}(i,j)\\
        &= (\widetilde{P}g)(i).
    \end{align*}
    
    For every eigenvalue $\lambda_i(\widetilde{P})$, let $f_i \in l^2(\overline{\pi})$ be an associated eigenfunction. With the fact that $G(Uf_i) = Uf_i$ since $Uf_i \in S$,
    \begin{align*}
        G Q_{\widetilde{P}} G (Uf_i) &= UU^* Q_{\widetilde{P}} Uf_i\\ 
        &= U\widetilde{P}f_i\\
        &= \lambda_i(\widetilde{P}) (Uf_i).
    \end{align*}
    Hence, every eigenvalue of $G Q_{\widetilde{P}} G (Uf_i)$ on $S$ is an eigenvalue of $\widetilde{P}$. On $S^\perp$, the eigenvalues must be $0$ since $G$ annihilates $S^\perp$. Hence, 
    $$\normalfont \mathrm{spec}(G Q_{\widetilde{P}} G)  = \mathrm{spec}(\widetilde{P}) \cup \{0\}.$$

\end{proof}

\begin{prop} \label{spectral sandwich = orbit}
    Let $P \in \mathcal{L}(\pi)$ be a sampler of $\mathcal{X}.$ For some non-trivial group action $\mathcal{G}$ and its orbits $(\mathcal{O}_i)_{i=1}^k$, where $k < n$, define the projection chain $\overline{P}$ as per \eqref{Pbar}. Then
    $$\normalfont \mathrm{spec}(GPG) = \mathrm{spec}(\overline{P}) \cup \{0\} .$$
\end{prop}

\begin{proof}
    Let $f \in \ell^2(\overline{\pi}).$ Then 
    \begin{align*}
        (U^*PUf)(i) &= \frac{1}{\pi(\mathcal{O}_i)} \sum_{x\in \mathcal{O}_i} \pi(x) (PUf)(x)\\
        &= \frac{1}{\pi(\mathcal{O}_i)} \sum_{x\in \mathcal{O}_i} \pi(x) \sum_{y\in \mathcal{X}} P(x,y)(Uf)(y)\\
        &= \frac{1}{\pi(\mathcal{O}_i)} \sum_{x\in \mathcal{O}_i} \pi(x) \sum_{j=1}^k \sum_{y \in \mathcal{O}_j} P(x,y) f(j)\\
        &= \sum_{j=1}^k \overline{P}(i,j) f(j),
    \end{align*}
    and hence $U^* PU = \overline{P}$.
    
    Suppose $h \in \ell^2(\pi)$ is an eigenfunction of $GPG$ with eigenvalue $\lambda \neq 0$. Then 
    $$\lambda Gh = G^2PGh = GPGh = \lambda h.$$
    Hence, $h \in S = \text{Im}(U)$ and so we can find $f \in \ell^2(\overline{\pi})$ such that $h = Uf.$

    It follows that 
    \begin{align*}
        \lambda f = U^* GPG Uf = U^* (UU^*) P (UU^*) Uf = U^*PUf = \overline{P}f,
    \end{align*}
    or equivalently, $\mathrm{spec}(GPG)\setminus \{0\}$ = $\mathrm{spec}(\overline{P}) \setminus \{0\}$. By similar argument, any eigenvalue $\lambda \neq 0$ corresponding to $\overline{P}$ must also be an eigenvalue of $GPG$.

    Since $\mathcal{G}$ admits $k < n$ orbits, $0$ must be an eigenvalue of $GPG$ as well. Hence, $\mathrm{spec}(GPG) = \mathrm{spec}(\overline{P}) \cup \{0\} .$
\end{proof}

Similarly, one can also look at $\overline{P}$ and $\overline{\Pi}$ to determine the KL-divergence of $GPG$ from $\Pi$.\\

\begin{prop} \label{KL orbit = sandwich}
    Given a group action $\mathcal{G}$ and its orbits $(\mathcal{O}_i)_{i=1}^k$, let $GPG$ be the Gibbs-orbit sampler associated with some sampler $P \in \mathcal{S}(\pi).$ Then 
    $$D_{KL}^\pi(GPG \| \Pi) = D_{KL}^{\overline{\pi}}\big(\overline{P} \| \overline{\Pi}\big),$$
    where $\overline{\Pi}$ is the matrix with each row equals to $\overline{\pi}$.
\end{prop}

\begin{proof}
    Using \eqref{GPGmatrix},

    \begin{align*}
        D_{KL}^\pi(GPG \| \Pi) &= \sum_{x,y \in \mathcal{X}} \pi(x) GPG(x,y) \log\bigg(\frac{GPG(x,y)}{\pi(y)}\bigg)\\
        &= \sum_{i,j=1}^k \sum_{\substack{x\in \mathcal{O}_i\\ y \in \mathcal{O}_j}} \pi(x) GPG(x,y) \log\bigg(\frac{1}{\pi(\mathcal{O}_i)\pi(\mathcal{O}_j)} \sum_{\substack{z \in \mathcal{O}(x)\\ w \in \mathcal{O}(y)}} \pi(z) P(z,w)\bigg)\\
        &= \sum_{i,j=1}^k \log\bigg(\frac{\overline{P}(i,j)}{\pi(\mathcal{O}_j)}\bigg) \sum_{\substack{x\in \mathcal{O}_i\\ y \in \mathcal{O}_j}}  \frac{\pi(x)\pi(y)}{\pi(\mathcal{O}_i)\pi(\mathcal{O}_j)}\sum_{\substack{z \in \mathcal{O}_i \\ w \in \mathcal{O}_j}} \pi(z)P(z,w)\\
        &= \sum_{i,j=1}^k \log\bigg(\frac{\overline{P}(i,j)}{\pi(\mathcal{O}_j)}\bigg) \overline{P}(i,j) \sum_{\substack{x\in \mathcal{O}_i\\ y \in \mathcal{O}_j}} \frac{\pi(x)\pi(y)}{\pi(\mathcal{O}_j)}\\
        &= \sum_{i,j=1}^k \pi(\mathcal{O}_i) \overline{P}(i,j) \log\bigg(\frac{\overline{P}(i,j)}{\pi(\mathcal{O}_j)}\bigg) \\
        &= D_{KL}^{\overline{\pi}}\big(\overline{P} \| \overline{\Pi}\big).
    \end{align*}
\end{proof}

Propositions \ref{spectral sandwich = orbit} and
\ref{KL orbit = sandwich} show that, for a fixed orbit partition, both the
spectrum and the KL divergence of $GPG$ are determined entirely by the
projection chain $\overline P$. Thus, candidate samplers may be designed on
the lower-dimensional orbit space and subsequently lifted to the original
state space using \eqref{liftingP}. We next illustrate this principle with a
particular choice of $\overline P$.

\begin{prop} \label{starorbit}
    Let $(\mathcal{O}_i)_{i=1}^k$ be the orbits given by a fixed group action $\mathcal{G}$, and suppose they are ordered $\pi(\mathcal{O}_1) \leq \dots \leq \pi(\mathcal{O}_k)$, with $\pi(\mathcal{O}_k) > 1/2.$ Then the sampler 
    $$\overline{P} = \begin{pmatrix}
        0 & 0 & \cdots & 0 & 1\\
        \vdots & \vdots & \ddots & \vdots & \vdots\\
        0 & 0 & \cdots & 0 & 1\\
        \frac{\pi(\mathcal{O}_1)}{\pi(\mathcal{O}_k)} & \frac{\pi(\mathcal{O}_2)}{\pi(\mathcal{O}_k)} & \cdots & \frac{\pi(\mathcal{O}_{k-1})}{\pi(\mathcal{O}_k)} & 2 - \frac{1}{\pi(\mathcal{O}_k)}
    \end{pmatrix}$$
    and its Gibbs sampler $GPG$ has absolute spectral gap $2-\pi(\mathcal{O}_k)^{-1}$. As $\pi(\mathcal{O}_k) \rightarrow 1$, we also have that $D_{KL}^{\overline{\pi}}(\overline{P} \| \overline{\Pi}) \rightarrow 0.$
\end{prop}

\begin{proof}
    Of the $k$ eigenvalues, $k-2$ of them will be 0 since $\mathrm{rank}(\overline{P}) = 2$. Then apart from the trivial eigenvalue 1, the last remaining eigenvalue would be $1 - \pi(\mathcal{O}_k)^{-1}$ with eigenvector $(1, \dots, 1, 1-\pi(\mathcal{O}_k)^{-1}).$

    The absolute spectral gap then follows from the fact that 
    $$\mathrm{spec}(\overline{P}) = \{1, 0, 1 - \pi(\mathcal{O}_k)^{-1}\},$$
    together with Proposition \ref{spectral sandwich = orbit}.

    The KL divergence of $\overline{P}$ from $\overline{\Pi}$ is 
    \begin{align*}
        D_{KL}^{\overline{\pi}}(\overline{P} \| \overline{\Pi}) &= \sum_{i,j=1}^k \overline{\pi}(i) \overline{P}(i,j) \log \bigg(\frac{\overline{P}(i,j)}{\overline{\pi}(j)}\bigg)\\
        &=  2 (1 - \overline{\pi}(\mathcal{O}_k)) \log \bigg( \frac{1}{\overline{\pi}(\mathcal{O}_k)}\bigg) + (2\overline{\pi}(\mathcal{O}_k)-1) \log \bigg(\frac{2\overline{\pi}(\mathcal{O}_k)-1}{\overline{\pi}(\mathcal{O}_k)^2}\bigg). 
    \end{align*}
    Hence, as $\overline{\pi}(\mathcal{O}_k) \rightarrow 1$, the expression goes to 0.
\end{proof}

Consider the family of samplers
\begin{align*}
    \mathcal{D} = \mathcal{D}(\mathcal{G},\pi) := \{P \in \mathcal{L}(\pi);~P(x,y) = 0 \,\, \mathrm{for \,\, all\,\,} x \in \mathcal{O}_i, y \in \mathcal{O}_j, i,j \in \llbracket k-1 \rrbracket\}.
\end{align*}

Every $P\in\mathcal D$ induces the same projection chain $\overline P$ from
Proposition~\ref{starorbit}. Consequently, all such lifts produce the same
value of $D_{KL}^\pi(GPG\|\Pi)$. The set $\mathcal D$ therefore describes a
family of state-space realisations of this particular orbit-space kernel.

Using \eqref{liftingP}, one such feasible $P$ is given by $Q_{\overline{P}} \in \mathcal{D}$ that lifts $\overline{P}$ back to the state space $\mathcal{X}$. Formally,

$$Q_{\overline{P}}(x,y) = \begin{cases}
    \frac{\pi(y)}{\pi(\mathcal{O}_k)}, & x \notin \mathcal{O}_k,\ y \in \mathcal{O}_k \ \mathrm{or}\ x \in \mathcal{O}_k,\ y \notin \mathcal{O}_k,\\[0.6ex]
    \frac{\pi(y) (2\pi(\mathcal{O}_k)-1)}{\pi(\mathcal{O}_k)^2}, & x,\ y \in \mathcal{O}_k,\\
    0, & \mathrm{otherwise}.
\end{cases}$$

\subsection{An example on the Curie-Weiss model}
We recall the mean-field Curie-Weiss model as described in Chapter 13 of \cite{Bovier_denHollander_2015}. The model is a high-dimensional system that has been widely studied in statistical mechanics and probability. 

Let the state space be $\mathcal{X} = \{-1, +1\}^d,$ for some positive even integer $d$. Then each configuration $x = (x_1, \dots, x_d)$ represents the spin orientation of $d$ interacting particles. The Hamiltonian of the model is given by 
$$H_d(x) = -\frac{1}{2d}\sum_{i,j = 1}^d x_i x_j - h \sum_{i =1}^d x_i,$$
with $h \in \mathbb{R}$ as the magnetic field. We shall assume $h=0$ for the rest of this subsection.

The Hamiltonian only depends on $x$ through its magnetisation  
$$m_d(x) = \frac{1}{d}\sum_{i=1}^d x_i.$$
That is,
$$H_d(x) = -\frac{d}{2} m^2_d(x).$$

The corresponding Gibbs distribution at inverse temperature $\beta>0$ is then, for $x \in \mathcal{X}$,
$$\pi_\beta(x) = \frac{1}{Z(\beta,d)} \exp(-\beta H_d(x)),$$
with $Z(\beta, d)$ as the normalising constant. 

Hence, $\pi_\beta(x)$ depends only on $m_d(x)$, and any pair $x,y \in \mathcal{X}$ with $m_d(x) = m_d(y)$ must have the same probability under $\pi_\beta$. Further, the model is invariant to the global flip spin $x \to -x$.

This motivates us to consider the partitions $(\mathcal{O}_i)_{i=0}^{d/2}$, where for $i \in \llbracket 0,d/2 \rrbracket$
$$\mathcal{O}_{i} = \bigg\{x \in \mathcal{X}: |m_d(x)| = \frac{2i}{d}\bigg\}.$$
Under each partition $\mathcal{O}_i$, all elements are uniformly distributed. Let
\[c_i= \begin{cases}
    1, & i=0,\\
    2, & i\geq1.
\end{cases}\]

Then
\[\pi_\beta(\mathcal O_i) = \frac{c_i}{Z(\beta,d)} \binom{d}{d/2-i} \exp\bigg(\frac{2\beta i^2}{d}\bigg).\]

Hence
\[\frac{\pi_\beta(\mathcal O_{i+1})} {\pi_\beta(\mathcal O_i)} =
\begin{cases}
    \dfrac{2d}{d+2}\exp\bigg(\dfrac{2\beta}{d}\bigg), & i=0,\\[1.4ex]
    \dfrac{\frac d2-i}{\frac d2+i+1} \exp\bigg(\dfrac{2\beta(2i+1)}{d}\bigg), & 1\leq i\leq d/2-1.
\end{cases}\]

Following which, a sufficient condition for monotonicity is to study the map $f: [0, d/2-1] \to \mathbb{R}$ defined by 
$$f(x) = \frac{d-2x}{d+2x+2}\exp\bigg(\frac{2\beta}{d}(2x+1)\bigg).$$
Take $g = \log f$, where
$$g(x) = \log(d-2x) - \log(d+2x+2) + \frac{2\beta}{d}(2x+1).$$

Its derivative 
\begin{align*}
    g'(x) &= -\frac{2}{d-2x} -\frac{2}{d+2x+2} + \frac{4\beta}{d}\\
    &= \frac{-4(d+1)}{(d-2x)(d+2x+2)} + \frac{4\beta}{d}
\end{align*}
is decreasing in $x$ on $[0,d/2-1]$. Hence, for $g'(x) > 0,$ it suffices for 
$$\frac{-4(d+1)}{4d} + \frac{4\beta}{d} \geq 0,$$
or equivalently 
$$\beta \geq \frac{d+1}{4}.$$

Under this condition, 
\begin{align*}
    \dfrac{\pi_\beta(\mathcal{O}_{i+1})}{\pi_\beta(\mathcal{O}_i)} = f(i) \geq f(0) = \frac{d}{d+2}\exp\bigg(\frac{2\beta}{d}\bigg) \geq \frac{d}{d+2}\bigg(1 + \frac{2\beta}{d}\bigg) \geq 1,
\end{align*}
so long as $\beta \geq 1$. Set $\beta^* = \max\{(d+1)/4, 1\}$. Then, at sufficiently large $\beta \geq \beta^*$, we then have $\pi_\beta(\mathcal{O}_0) \leq \pi_\beta(\mathcal{O}_1) \leq \dots \leq \pi_\beta(\mathcal{O}_{d/2}).$ 

Now consider the projection chain induced by $(\mathcal{O}_i)_{i=0}^{d/2}$. Let $\overline{\pi}_\beta = (\pi_\beta(\mathcal{O}_0), \dots , \pi_\beta(\mathcal{O}_{d/2}))$ and suppose we seek $\overline{G}$, the best Gibbs kernel on $\overline{\pi}_\beta$ in terms of KL divergence to $\overline{\Pi}_\beta$. 

Proposition \ref{starorbit} proposes the following orbit $(\mathcal{B}_r)_{r=0}^k$, with $k \in \llbracket d/2 \rrbracket$:
$$\mathcal{B}_r = \begin{cases}
    \mathcal{O}_{r}, & \mathrm{if}\ r \in \llbracket 0,k-1 \rrbracket,\\
    \mathcal{O}_k\cup \dots \cup \mathcal{O}_{d/2}, & \mathrm{if}\ r = k.
\end{cases}$$

\begin{figure}[ht]
  \centering
  \begin{tikzpicture}
  \begin{axis}[
    axis lines=middle,
    xlabel={$m$}, ylabel={$H_d(m)$},
    xmin=-1.1, xmax=1.2,
    ymin=-0.05, ymax=1.2,
    xtick=\empty, ytick=\empty,
    axis line style={thick},
    samples=200, domain=-1:1,
    width=10cm, height=6cm,
    clip=false
  ]

    % --- Potential curve ---
    \addplot[thick,blue!70!black] {1 - x^2};

    % --- Dashed intermediate B_r levels ---
    \addplot[dashed,gray,domain=-0.98:0.98] {0.9};
    \node[anchor=west] at (axis cs:1.03,0.9) {$\mathcal{B}_1$};

    \addplot[dashed,gray,domain=-0.98:0.98] {0.7};
    \node[anchor=west] at (axis cs:1.03,0.7) {$\mathcal{B}_2$};

    \node[anchor=west] at (axis cs:1.03,0.575) {$\vdots$};

    \addplot[dashed,gray,domain=-0.98:0.98] {0.4};
    \node[anchor=west] at (axis cs:1.03,0.4) {$\mathcal{B}_{k-1}$};

    % --- Single adjointed shaded box for B_k ---
    \addplot [draw=gray!70, fill=gray!25, opacity=0.6] 
      coordinates {(-1,0.0) (-1,0.3) (1,0.3) (1,0.0)} -- cycle;
    \node[anchor=west] at (axis cs:1.03,0.2) {$\mathcal{B}_k$};

    % --- Custom ±1 labels ---
    \node[below] at (axis cs:-1,0) {$-1$};
    \node[below] at (axis cs: 1,0) {$+1$};

  \end{axis}
\end{tikzpicture}
  \caption{Plot of $H_d(m)$ against different magnetisation levels and the orbits $\mathcal{B}_r$.}
  \label{fig:curie-weiss}
\end{figure}

As $\beta \rightarrow \infty$, the mass of $\pi_\beta$ increasingly concentrates about the modes $+1$ and $-1$. This implies $\pi_\beta(\mathcal{B}_k) \rightarrow 1$, and so, the sampler $\overline{P}$ described in Proposition \ref{starorbit} is a suitable candidate for sampling over $(\mathcal{B}_r)_{r=0}^k.$ 

In fact, as $\beta \rightarrow \infty$, one can take $k = d/2$, that is, to use the original orbits $(\mathcal{O}_i)_{i=0}^{d/2}$, since the bulk of the mass would be concentrated at $\pi_\beta(\mathcal{O}_{d/2}).$

After which, one can formulate the sampler $Q_{\overline{P}}$ similar to \eqref{liftingP} as
\begin{equation} \label{starproj}
    Q_{\overline{P}}(x,y) = \begin{cases}
    \frac{\pi_\beta(y)}{\pi_\beta(\mathcal{B}_k)}, & x \notin \mathcal{B}_k,\ y \in \mathcal{B}_k \ \mathrm{or}\ x \in \mathcal{B}_k,\ y \notin \mathcal{B}_k,\\[0.6ex]
    \frac{\pi_\beta(y) (2\pi_\beta(\mathcal{B}_k)-1)}{\pi_\beta(\mathcal{B}_k)^2}, & x,\ y \in \mathcal{B}_k,\\
    0, & \mathrm{otherwise}.
\end{cases}
\end{equation}

In practice, this is how one could implement $Q_{\overline{P}}$.

If the sampler is at state $x \notin \mathcal{B}_k$, the construction of $\overline{P}$ guarantees that the next state $y$ will be in $\mathcal{B}_k$. Then 
\begin{enumerate}
    \item Draw an orbit index $i \in \llbracket k, d/2 \rrbracket$ with $\mathbb{P}(i) = \pi_\beta(\mathcal{O}_i)/ \pi_\beta(\mathcal{B}_k).$
    \item Draw $y$ uniformly within the orbit $\mathcal{O}_i$.
\end{enumerate}

If the sampler is at state $x \in \mathcal{B}_k$, one of two cases could happen. With probability $2-\pi_\beta(\mathcal{B}_k)^{-1}$, the next state will be within $\mathcal{B}_k$. Then 
\begin{enumerate}
    \item Draw an orbit index $i \in \llbracket k, d/2 \rrbracket$ with $\mathbb{P}(i) = \pi_\beta(\mathcal{O}_i)/ \pi_\beta(\mathcal{B}_k).$
    \item Draw $y$ uniformly within the orbit $\mathcal{O}_i$.
\end{enumerate}
Else, the next jump will be to some $y \in \mathcal{O}_i$ for $i \in \llbracket k-1 \rrbracket$. Then 
\begin{enumerate}
    \item Draw an orbit index $i \in \llbracket k-1 \rrbracket$ with $\mathbb{P}(i) = \pi_\beta(\mathcal{O}_i)/ (1-\pi_\beta(\mathcal{B}_k)).$
    \item Draw $y$ uniformly from $\mathcal{O}_i$.
\end{enumerate}

One way to draw $y$ uniformly from an orbit $\mathcal{O}_i$ is to utilise the Fisher-Yates (or Knuth) shuffle algorithm as described in Chapter 3.4.2 of \cite{Knuth_1997}. One can use the algorithm to sample a permutation $y \in \mathcal{X}$ with $d/2 + i$ number of $+1$'s. Then, sample a sign from $\{\pm1\}$ to return either $+y$ or $-y$.

A key distinction from previous work by \cite{Choi_Permutation} is that only equi-probability permutations were considered. That is, within each orbit, every state must have exactly the same probability with respect to the target $\pi$. Now however, we are able to group states with similar but not necessarily equal probabilities together into an orbit. The resulting Gibbs sampler $GPG$ would then improve mixing over any original sampler $P$.

Finally, we discuss the performance of our sampler in comparison to the usual Metropolis-Hastings algorithm on $\mathcal{X}$.

Recall from \cite{Levin_Peres_Wilmer_2017}, we define the total variation distance between any $\mu, \nu \in \mathcal{P}(\mathcal{X}),$
$$\| \mu - \nu \|_{TV} := \frac{1}{2} \sum_{x\in \mathcal{X}} |\mu(x)-\nu(x)|,$$
and the worst-case total variation mixing time of the Markov chain associated with P, for some $\epsilon > 0$, is
$$t_{mix}(P,\varepsilon) := \inf\bigg\{ n \in \mathbb{N}:\ \max_{x\in \mathcal{X}} \|P^n(x,\cdot) - \pi \|_{TV} < \varepsilon\bigg\}.$$

Define the relaxation time of a reversible Markov chain with absolute spectral gap $\lambda$ as
$$t_{rel} := \frac{1}{\lambda}.$$

Then by Theorem 12.4 of \cite{Levin_Peres_Wilmer_2017}, 
$$t_{mix}(Q_{\overline{P}}, \varepsilon) \leq t_{rel}\big(\overline{P}\big)\log \bigg(\frac{1}{\varepsilon \pi_{\min}}\bigg),$$
where $\pi_{\min} = \min\{\pi_\beta(x): x \in \mathcal{X}\}.$

Notice that 
\begin{align*}
     Z(\beta,d) &= \binom{d}{d/2} + 2\sum_{i=1}^{d/2} \binom{d}{d/2-i} \exp\bigg(\frac{2\beta i^2}{d}\bigg)\\
    &\leq \exp\bigg(\frac{d}{2}\beta\bigg) \sum_{i=0}^d \binom{d}{i}\\
    &= 2^d\exp\bigg(\frac{d}{2}\beta\bigg).
\end{align*}
Since $\pi_{\min} = Z(\beta,d)^{-1}$, $\lambda(Q_{\overline{P}}) = 2-\pi_\beta(\mathcal{B}_k)^{-1}$, and by writing $\pi_\beta(\mathcal{B}_k) = \frac{1}{2} + \delta$ for $\delta > 0$, we see that for $\beta \geq \beta^*$,
$$t_{mix}(Q_{\overline{P}}, \varepsilon) \leq \frac{\pi_\beta(\mathcal{B}_k)}{2\pi_\beta(\mathcal{B}_k) - 1} \bigg(\frac{d\beta}{2} + d\log(2) - \log(\varepsilon)\bigg) \leq \frac{1}{2\delta} \bigg(\frac{d\beta}{2} + d\log(2) - \log(\varepsilon)\bigg).$$

This implies that the mixing time of $Q_{\overline{P}}$ is at most polynomial in $\beta$, $d$ and $1/\delta$. 

In contrast, a classical sampler $P$ in this context is the Glauber dynamics that targets $\pi_\beta$. That is, at each iteration a coordinate is chosen uniformly-at-random out of the $d$ coordinates and is flipped to the opposite spin. This proposal configuration is then subjected to a Metropolis-Hastings filter that targets $\pi_\beta$. For such $P$, by Theorem $12.5$ in \cite{Levin_Peres_Wilmer_2017}, we note that
\begin{align*}
    t_{mix}(P, \varepsilon) &\geq \left(\dfrac{1}{1-\lambda_2(P)} - 1\right) \log \left(\dfrac{1}{2 \varepsilon} \right) \\
    &\geq \left(\dfrac{e^{\beta d}}{4^d} - 1\right) \log \left(\dfrac{1}{2 \varepsilon} \right),
\end{align*}
where the last inequality follows from Lemma $2.3$ in \cite{Holley_Stroock_1988}. This implies that the mixing time of $P$ is at least exponential in $\beta$ and $d$ for $\beta \geq \beta^*$. 

We summarize the above discussions into the following proposition.\\

\begin{prop}
For the mean-field Curie--Weiss model, fix $\beta \geq \beta^* = \max\{(d+1)/4,\,1\}$. Let $Q_{\overline{P}}$ be defined as in \eqref{starproj}, and let $P$ denote the single-site Glauber dynamics targeting $\pi_\beta$. 
Then the worst-case total variation mixing times satisfy
$$
t_{\mathrm{mix}}(Q_{\overline{P}}, \varepsilon)\leq \frac{1}{2\delta} \bigg(\frac{d\beta}{2} + d\log(2) - \log(\varepsilon)\bigg),
$$
where $\pi_\beta(\mathcal{B}_k) = \frac{1}{2} + \delta$, while
$$t_{mix}(P, \varepsilon) \geq \left(\dfrac{e^{\beta d}}{4^d} - 1\right) \log \left(\dfrac{1}{2 \varepsilon} \right).$$
In particular, $t_{\mathrm{mix}}(Q_{\overline{P}},\varepsilon)$ is at most polynomial in $d$, $\beta$ and $1/\delta$, whereas $t_{mix}(P, \varepsilon)$ is at least exponential in $d$ and $\beta$.
\end{prop}

The partitioning of the Curie–Weiss model by magnetisation can be viewed as a discrete analogue of the energy rings used in the equi-energy sampler of \cite{Kou_2006}. In their framework, the Gibbs measure is decomposed into level sets of the Hamiltonian, and transitions are designed to exchange information between states of comparable energy that are otherwise separated by steep energy barriers. This is closely analogous to our construction of group orbits by magnetisation, where the state space is stratified into groups, each corresponding to an exact energy level, to promote efficient mixing across modes of similar potential.

\section{Optimal Gibbs projection in KL divergence and exact-sampling conditions}\label{section:bestG}

We now optimise the Gibbs orbit kernel $G$ itself. Specifically, among orbit partitions with a fixed number $k$ of blocks, we identify the Gibbs projection minimizing $D_{KL}^\pi(G\|\Pi)$.

\begin{prop} \label{bestGbyKL}
    Suppose $\pi$ is ordered in non-decreasing order, that is, $\pi(1) \leq \pi(2) \leq \dots \leq \pi(n),$ and for $k \in \llbracket n \rrbracket$, let $(\mathcal{O}_i)_{i=1}^k$ be the partition 
    \begin{equation} \label{bestorbitKL}
        \mathcal{O}_i = \begin{cases}
        \{ i \}, & \mathrm{if}\ 1 \leq i \leq k-1,\\
        \{ k, k+1, \dots, n\}, & \mathrm{if}\ i = k.
    \end{cases}
    \end{equation}

    Then the Gibbs kernel $G_\mathcal{O}$ defined by $(\mathcal{O}_i)_{i=1}^k$ is the minimiser of 
    $D_{KL}^\pi(G \| \Pi)$ among all other Gibbs kernel with $k$ orbits. That is, for any other orbit $(\mathcal{C}_i)_{i=1}^k$, 
    $$D_{KL}^\pi(G_\mathcal{O} \| \Pi) \leq D_{KL}^\pi(G_\mathcal{C} \| \Pi).$$
\end{prop}

\begin{proof}
    If $k=n$, then $(\mathcal{O}_i)$ is the only permissible partition and the claim holds trivially.

    Fix $k \in \llbracket n \rrbracket$. Let
    $$H(\pi) := -\sum_{x\in\mathcal X} \pi(x)\log\pi(x)$$
    denote the Shannon entropy of $\pi$. For a given partition $(\mathcal{O}_i)_{i=1}^k$, recall $\overline{\pi} = (\pi(\mathcal{O}_1),\ldots,\pi(\mathcal{O}_k)).$
    Then we define
    $$H(\overline{\pi}):= -\sum_{i=1}^k \pi(\mathcal{O}_i)\log \pi(\mathcal{O}_i)$$
    be the entropy of the corresponding block masses.
    
    By Proposition \ref{prop:PythidGPG},
    $$D_{KL}^\pi(I \| \Pi) = D_{KL}^\pi(I \| G) + D_{KL}^\pi(G \| \Pi),$$
    and note that 
    \begin{align*}
        D_{KL}^\pi(I \| G) &= \sum_{x\in \mathcal{X}} \pi(x) \log \frac{1}{\pi(x)/\pi(\mathcal{O}(x))}\\
        &= H(\pi) + \sum_{i=1}^k \pi(\mathcal{O}_i) \log (\pi(\mathcal{O}_i))\\
        &= H(\pi) - H(\overline\pi).
    \end{align*}
    Hence, for fixed $\pi$, minimising $D_{KL}^\pi(G \| \Pi)$ is equivalent to minimising $H(\overline\pi)$ over all partitions with $k$ blocks.

    Let $g(t)=t\log t$ . For any two blocks $\mathcal{C}_i, \mathcal{C}_j$ with total mass $S = \pi(\mathcal{C}_i)+\pi(\mathcal{C}_j)$, define $h(t)=g(t)+g(S-t)$. Then $h'' > 0$ on $(0,S)$, and so $h$ is strictly convex on the same interval and achieves its maximum at the endpoints.

    Now let $(\mathcal{C}_i)_{i=1}^k$ be any partition that differs from $(\mathcal{O}_i)_{i=1}^k$. Because $k<n$, there exists at least one non-singleton block, denoted by $\mathcal{C}_i$.

    Suppose there exists another non-singleton block $\mathcal{C}_j$, and we let $x_m$ be the element within $\mathcal{C}_i \cup \mathcal{C}_j$ with the smallest probability. Since the two blocks are non-singletons, their masses lie strictly between $\pi(x_m)$ and $S-\pi(x_m).$ 
    
    By the strict convexity of $h,$
    $$g(\pi(\mathcal{C}_i)) + g(\pi(\mathcal{C}_j)) = h(\pi(\mathcal{C}_i)) < h(\pi(x_m)) = g(\pi(x_m))+g(\pi(\mathcal{C}_i \cup \mathcal{C}_j) \setminus \{x_m\}).$$
    
    Thus replacing the pair $(\mathcal{C}_i,\mathcal{C}_j)$
    by a new pair consisting of the singleton $\{x_m\}$ and
    the merged remainder
    $(\mathcal{C}_i \setminus \{x_m\}) \cup \mathcal{C}_j$
    strictly decreases $H(\overline\pi)$.

    Iterating this push-out operation would produce a partition with exactly one non-singleton block and $k-1$ singletons.

    Suppose the partition now has a single non-singleton $\mathcal{C}_i$, and all remaining blocks are singletons. If every singleton $\mathcal{C}_j$ satisfies, 
    $$\pi(\mathcal{C}_j) \leq \min_{x \in \mathcal{C}_i} \pi(x),$$
    then $(\mathcal{C}_i)_{i=1}^k = (\mathcal{O}_i)_{i=1}^k$. 
    
    Otherwise choose a singleton $\mathcal{C}_j = \{y\}$ such that $\pi(y) > \pi(x_m) = \min_{x \in \mathcal{C}_i} \pi(x).$ 

    By the same convexity argument as before, we can show that replacing $\mathcal{C}_i$ and $\{y\}$ by the pair $\{x_m\}$ and $(\mathcal{C}_i \setminus \{x_m\}) \cup \mathcal{C}_j$ will again strictly decreases the entropy. 

    Repeating such swaps eventually yields a configuration in which all singletons correspond to the 
    $k-1$ smallest atoms, or equivalently, the orbits described as $(\mathcal{O}_i)_{i=1}^k$.
\end{proof}

Having identified the KL-optimal Gibbs projection $G$, we next ask when
augmentation by this $G$ yields exact sampling. Equivalently, we characterize
the base kernels $P$ for which $GPG=\Pi$.

\begin{prop} \label{bestPgivenG}
    Suppose $\pi(1) \leq \dots \leq \pi(n)$. Let $G$ be the Gibbs kernel constructed by $(\mathcal{O}_i)_{i=1}^k$ defined in \eqref{bestorbitKL}. Then $D_{KL}^\pi(GPG \| \Pi) = 0$ if and only if $P$ satisfies the following conditions:
    \begin{enumerate}
        \item $P(x,y) = \pi(y)\ \mathrm{for}\ x,y \in \llbracket k-1 \rrbracket.$

        \item $\sum_{w \in \mathcal{O}_k} P(x,w) = \pi(\mathcal{O}_k)\ \mathrm{for}\ x \in \llbracket k-1 \rrbracket.$

        \item $\sum_{z\in \mathcal{O}_k} \pi(z) P(z,y) = \pi(\mathcal{O}_k) \pi(y)\ \mathrm{for}\ y \in \llbracket k-1 \rrbracket.$

        \item $\sum_{z,w \in \mathcal{O}_k} \pi(z) P(z,w) = (\pi(\mathcal{O}_k))^2.$
    \end{enumerate}
     Equivalently, this implies that $GPG = \Pi$.
\end{prop}

\begin{proof}
    Recall that in \eqref{GPGmatrix},
    $$GPG(x,y) = \frac{\pi(y)}{\pi(\mathcal{O}(x)) \pi(\mathcal{O}(y))} \sum_{\substack{z \in \mathcal{O}(x) \\ w \in \mathcal{O}(y)}} \pi(z) P(z,w).$$
    Now consider the 4 cases:

    For $x,y \in \llbracket k-1 \rrbracket$, $x$ and $y$ are in their respective singleton orbits. Then
    $$GPG(x,y) = \frac{\pi(x)}{\pi(x)\pi(y)}\pi(x)P(x,y) = P(x,y).$$

    If $x \in \llbracket k-1 \rrbracket$ and $y \in \mathcal{O}_k$,
    $$GPG(x,y) = \frac{\pi(y)}{\pi(x)\pi(\mathcal{O}_k)} \sum_{w\in \mathcal{O}_k}\pi(x) P(x,w).$$

    If $x \in \mathcal{O}_k$ and $y \in \llbracket k-1 \rrbracket,$
    $$GPG(x,y) = \frac{1}{\pi(\mathcal{O}_k)} \sum_{z \in \mathcal{O}_k} \pi(x)P(x,y).$$

    Lastly, if both $x,y \in \mathcal{O}_k$,
    $$GPG(x,y) = \frac{\pi(y)}{\pi(\mathcal{O}_k)^2}\sum_{z,w \in \mathcal{O}_k} \pi(z)P(z,w).$$

    The four conditions listed then follows from the fact that for $GPG = \Pi$, $GPG(x,y) = \pi(y)$ for all $x,y \in \mathcal{X}$. 
\end{proof}

\begin{remark}
    Note that the family of $P \in \mathcal{S}(\pi)$ described by Proposition \ref{bestPgivenG} does not only contain $\Pi$. We describe a class of such $P \neq \Pi$ that satisfies $GPG = \Pi$, where $G$ is constructed as per Proposition \ref{bestGbyKL}.

    Define $P$ such that for $y \in \llbracket k-1 \rrbracket,$ $P(x,y) = \pi(y)$, and for any 
     $x \in \llbracket k-1 \rrbracket$, the entries starting from column $k$ to $n$ can be arbitrary so long as they add up to $\pi(\mathcal{O}_k).$
    For $x,y \in \mathcal{O}_k$, we also set
    $$P(x,y) = \frac{1}{\pi(\mathcal{O}_k)}\bigg( \pi(y) - \sum_{x=1}^{k-1} \pi(x)P(x,y)\bigg).$$
\end{remark}

A concrete example is as follows:
Let $\pi = (0.05, 0.1, 0.2, 0.25, 0.4)$ and suppose $\mathcal{G}$ has orbits $\{1\},\ \{2\},\ \{3,4,5\}.$ Then 
$$P(x,y) = \begin{pmatrix}
    0.05 & 0.1 & 0 & 0.35 & 0.50 \\
    0.05 & 0.1 & 0.6 & 0.25 & 0 \\ 
    0.05 & 0.1 & 14/85 & 83/340 & 15/34 \\ 
    0.05 & 0.1 & 14/85 & 83/340 & 15/34 \\
    0.05 & 0.1 & 14/85 & 83/340 & 15/34
    \end{pmatrix} \neq \Pi.$$

\section{Alternating group actions} \label{section:altproj}
In previous sections, we have shown that the augmented kernel $GPG$ always performs no worse than $P$ in terms of absolute spectral gap, and that it does not increase the asymptotic variance of $G$-invariant observables. This motivates the concept of alternating group actions, where we consider several group actions and repeated augmentations. 

\subsection{Alternating projections on $k$ groups}
Let $\mathcal{G}_1, \dots \mathcal{G}_k$ be $k$ different groups that would act on $\mathcal{X}$, with their respective Gibbs kernel $G_1, \dots, G_k$. Then each $G_i$ is an orthogonal projection onto the subspace $S_i$, defined as\
\begin{equation} \label{S_i}
    S_i = \{f \in \ell^2(\pi)\ |\ f(x) = f(y) \text{ if } x,y \text{ are in the same orbits under } \mathcal{G}_i\}.
\end{equation}
These subspaces are all of finite dimensions, and are hence closed. 

Using known results in the literature of alternating projections, which \cite{ginat2018methodalternatingprojections} gives an extensive overview, the projection $G_\infty$ that satisfies
$$\lim_{n \rightarrow \infty}\|(G_1 G_2\cdots G_k)^n - G_\infty\|_{\ell^2(\pi) \rightarrow \ell^2(\pi)} = 0$$
exists, and is the projection onto the closed subspace $S = \bigcap_{i=1}^k S_i \subseteq \ell^2(\pi).$ One may understand $G_\infty$ to be the limiting projection of $(G_1 G_2\cdots G_k)^n$ in the operator norm sense.  

For two closed subspaces $S_1, S_2$, the cosine as defined by \cite{deutsch2001best}, Definition 9.4, is
\begin{align} \label{cosine}
    c(S_1, S_2) &:= \sup\{ \langle f,h \rangle_\pi\ |\ f\in S_1 \cap S^\perp,\ h \in S_2 \cap S^\perp,\ \|f\|_\pi, \|h\|_\pi \leq 1\}\\
    &= \|G_1G_2 - G_{S_1 \cap S_2}\|_{\ell^2(\pi) \rightarrow \ell^2(\pi)} \notag,
\end{align}
where $G_{S_1 \cap S_2}$ is the projection onto the intersection $S_1 \cap S_2$.

Then the rate of convergence, for the case where $k=2$ is given in \cite{deutsch2001best}, Definition 9.8 by
\begin{equation} \label{altconvergence2}
    \|(G_1G_2)^n - G_\infty\|_{{\ell^2(\pi) \rightarrow \ell^2(\pi)}} = c(S_1, S_2)^{2n-1}.
\end{equation}

For any arbitrary $k$, we can generalise the concept of cosine by
\begin{equation} \label{kcosine}
    c_i = c(S_i, \cap_{j=i+1}^k S_j) \text{ and } c := \bigg[1-\prod_{i=1}^{k-1} (1-c_i^2)\bigg]^{1/2},
\end{equation}
and the rate of convergence is instead given by
\begin{equation} \label{altconverencek}
    \|(G_1\cdots G_k)^n - G_\infty\|_{\ell^2(\pi) \rightarrow \ell^2(\pi)} \leq c^n.
\end{equation}

For ease of notation, let $P \in \mathcal{L}(\pi)$ and we set
\begin{align*}
    K_n &:= (G_1\cdots G_k)^nP(G_k \cdots G_1)^n,\ K_\infty := G_\infty P G_\infty, \text{ and } T := G_1\cdots G_k.
\end{align*}

\begin{prop} \label{rhoconvergencek}
    For any $k$ Gibbs kernels $G_1, G_2, \dots, G_k$, and its limiting projection $G_\infty$, 
    $$0 \leq\ \rho((G_1G_2 \cdots G_k)^n P (G_k \cdots G_2G_1)^n) - \rho(G_\infty PG_\infty) \leq 2c^n \rho(P).$$
\end{prop}

\begin{proof}
    Consider 
    $$K_n - K_\infty = (T^n - G_\infty)P(T^n)^* + G_\infty P((T^n)^* - G_\infty).$$
    Then since the operator norm is invariant under adjoint (see \cite{Rudin1991} Theorem 4.10), 
    $$
    \|(T^n)^* - G_\infty\|_{\ell^2_0(\pi) \rightarrow \ell^2_0(\pi)} = \|(T^n - G_\infty)^*\|_{\ell^2_0(\pi) \rightarrow \ell^2_0(\pi)} = \|T^n - G_\infty\|_{\ell^2_0(\pi) \rightarrow \ell^2_0(\pi)}.
    $$
    By the subadditivity and submultiplicativity properties of the operator norm, 
    \begin{align*}
        \|K_n - K_\infty\|_{\ell^2_0(\pi) \rightarrow \ell^2_0(\pi)} &\leq \rho(P) \big( \|T^n - G_\infty\|_{\ell^2_0(\pi) \rightarrow \ell^2_0(\pi)} + \|(T^n)^* - G_\infty\|_{\ell^2_0(\pi) \rightarrow \ell^2_0(\pi)}\big)\\
        &\leq 2\rho(P)\ c^n,
    \end{align*}
    with $c$ given by \eqref{kcosine}.
    
    Furthermore, by Proposition \ref{SLEMGPG}, one has that $\rho(K_n)$ decreases monotonically towards $\rho(K_\infty)$. Hence,
    $$0 \leq\ \rho(K_n) - \rho(K_\infty) \leq \|K_n - K_\infty\|_{\ell^2_0(\pi) \rightarrow \ell^2_0(\pi)} \leq 2\rho(P)\ c^n.$$
\end{proof}

\begin{prop} \label{varconvergencek}
    Under the same setting as Proposition~\ref{rhoconvergencek}, suppose that $P$ is ergodic. 
    
    Then, for every $f\in\ell^2_0(\pi)$,
    \begin{equation} \label{eq:variance-alternating-c}
        \left|v(f,K_n)-v(f,K_\infty)\right| \leq \frac{4\rho(P)}{(1-\rho(P))^2}\, c^n\|f\|_\pi^2.
    \end{equation}
    
    If, in addition, $f\in S\cap\ell^2_0(\pi)$, then
    \[v(f,K_{n+1})\leq v(f,K_n), \qquad n\geq1,\]
    
    and hence
    $$0 \leq v(f,K_n)-v(f,K_\infty) \leq \frac{4\rho(P)}{(1-\rho(P))^2}\, c^n\|f\|_\pi^2.$$
\end{prop}

\begin{proof}
    In the proof of Proposition \ref{rhoconvergencek}, we established that 
    $$\|K_n - K_\infty\|_{\ell^2_0(\pi) \rightarrow \ell^2_0(\pi)} \leq 2\rho(P)\ c^n.$$
    
    Because $\|K_n\|,\|K_\infty\|\leq \rho(P)<1$, both resolvents are well-defined on $\ell^2_0(\pi)$ and satisfy
    $$\|(I-K_n)^{-1}\|_{\ell^2_0(\pi) \rightarrow \ell^2_0(\pi)}, \; \|(I-K_\infty)^{-1}\|_{\ell^2_0(\pi) \rightarrow \ell^2_0(\pi)} \leq \frac{1}{1-\rho(P)}.$$
    
    By the resolvent identity,
    $$(I-K_n)^{-1}-(I-K_\infty)^{-1} = (I-K_n)^{-1} (K_n-K_\infty) (I-K_\infty)^{-1}.$$
    
    Consequently,
    \begin{align*}
        \left|v(f,K_n)-v(f,K_\infty)\right| &= 2\left| \left\langle f, \left[(I-K_n)^{-1}-(I-K_\infty)^{-1} \right]f \right\rangle_\pi\right|\\
        &\leq 2\, \|(I-K_n)^{-1}\|_{\ell^2_0(\pi) \rightarrow \ell^2_0(\pi)}\ \|K_n-K_\infty\|_{\ell^2_0(\pi) \rightarrow \ell^2_0(\pi)}\\
        &\qquad\times \|(I-K_\infty)^{-1}\|_{\ell^2_0(\pi) \rightarrow \ell^2_0(\pi)}\ \|f\|_\pi^2\\
        &\leq \frac{4\rho(P)}{(1-\rho(P))^2}\ c^n \|f\|_\pi^2.
    \end{align*}
    
    It remains to prove the monotonicity assertion for $f\in S\cap\ell_0^2(\pi)$. Since
    \[K_{n+1} = G_1\cdots G_k K_n G_k\cdots G_1,\]
    
    define $L_0:=K_n$ and, recursively,
    \[L_j = G_{k-j+1}L_{j-1}G_{k-j+1}, \qquad j=1,\ldots,k.\]
    
    Then $L_k=K_{n+1}$.
    
    Since $f\in S=\bigcap_{i=1}^kS_i$, we have $G_if=f$ for every $i=1,\ldots,k$. Hence, by repeated application of Proposition~\ref{varGPG/P},
    \[v(f,L_j)\leq v(f,L_{j-1}), \qquad j=1,\ldots,k.\]
    
    Therefore, $v(f,K_{n+1})\leq v(f,K_n).$ Moreover, \eqref{eq:variance-alternating-c} implies that $v(f,K_n)\to v(f,K_\infty).$
    
    Since the sequence is non-increasing, it follows that $0\leq v(f,K_n)-v(f,K_\infty).$
\end{proof}

As a corollary, we present the case $k=2$, for which $c = c(S_1, S_2)$.

\begin{corollary} \label{convergencek=2}
    Under the assumption of Propositions \ref{rhoconvergencek} and \ref{varconvergencek}, for $k=2$, 
    $$0 \leq \rho(K_n)-\rho(K_\infty) \leq 2\rho(P)\,c(S_1,S_2)^{2n-1}.$$
    
    Moreover, for every $f\in\ell_0^2(\pi)$,
    $$\left|v(f,K_n)-v(f,K_\infty)\right| \leq \frac{4\rho(P)}{(1-\rho(P))^2}\, c(S_1,S_2)^{2n-1}\|f\|_\pi^2.$$
    
    If $f\in S_1\cap S_2\cap\ell_0^2(\pi)$, then $v(f,K_n)\downarrow v(f,K_\infty)$
    and 
    $$0 \leq v(f,K_n)-v(f,K_\infty) \leq \frac{4\rho(P)}{(1-\rho(P))^2}\, c(S_1,S_2)^{2n-1}\|f\|_\pi^2.$$
\end{corollary}

\begin{proof}
    For $k=2$, the alternating-projection identity \eqref{altconvergence2} gives
    $$\|K_n - K_\infty\|_{\ell^2_0(\pi) \rightarrow \ell^2_0(\pi)} \leq 2\rho(P)\ c(S_1, S_2)^{2n-1}.$$
    
    The results then follow from the same proofs as Proposition \ref{rhoconvergencek} and \ref{varconvergencek}, with the new bound of $\|K_n - K_\infty\|_{\ell^2_0(\pi) \rightarrow \ell^2_0(\pi)}$.
\end{proof}

\subsection{Practical implementation of alternating projections}
While alternating projections can improve the mixing of a sampler, it is typically computationally infeasible to perform a large number of iterations. Hence, rather than taking iterating products of $G_1\cdots G_k$, we aim to identify the limiting projection $G_\infty$ by characterising the subspace $S = \bigcap_{i=1}^k S_i$ instead. 

Let $\mathcal{O}_i(x)$ denote the orbit of $x$ under $\mathcal{G}_i$.
We define an equivalence relation $\sim$ on $\mathcal{X}$ as the transitive
closure of belonging to the same $\mathcal{G}_i$-orbit. More precisely,
for $x,y\in\mathcal{X}$, we write $x\sim y$ if there exist
$x=x_0,x_1,\dots,x_r=y$ and indices
$i_1,\dots,i_r\in\llbracket k\rrbracket$ such that
\[x_\ell\in\mathcal{O}_{i_\ell}(x_{\ell-1}),\qquad \ell=1,\dots,r.\]

\begin{prop} \label{G_infty}
    The limiting projection $G_\infty$ projects onto the subspace defined by 
    $$S = \bigcap_{i=1}^k S_i = \{f \in \ell^2(\pi): f \text{ is constant on the equivalence classes of } \sim \}.$$
\end{prop}

\begin{proof}
    Take any $f \in S$, and let $x \sim y$. By definition, there exist
    $x=x_0,x_1,\dots,x_r=y$ and indices
    $i_1,\dots,i_r\in\llbracket k\rrbracket$ such that
    $x_\ell\in\mathcal{O}_{i_\ell}(x_{\ell-1})$ for each
    $\ell=1,\dots,r$. Since $f\in S_{i_\ell}$, it is constant on every
    $\mathcal{G}_{i_\ell}$-orbit, and hence $f(x_{\ell-1})=f(x_\ell)$
    for each $\ell$. Therefore $f(x)=f(y)$, so $f$ is constant on every
    equivalence class of $\sim$.

    Conversely, if $f$ is constant on every equivalence class of $\sim$,
    then any two points in the same $\mathcal{G}_i$-orbit are equivalent
    under $\sim$. Hence $f$ is constant on every orbit defined by
    $\mathcal{G}_i$, so $f\in S_i$ for all $i$, and therefore $f\in S$.
\end{proof}

The implication of Proposition \ref{G_infty} is that one can determine $G_\infty$ simply by determining the equivalence classes of $\mathcal{X}$ under $\sim$. 

A simple way to construct the equivalence class is as follows: Start from any $x \in \mathcal{X}$ and run through all orbits $O_i(x)$, adding every element within to the same class as $x$. Reiterate this procedure till all elements of $\mathcal{X}$ has been accounted for. 

With this, one can construct $G_\infty$ as per \eqref{matrixG}, taking the equivalence class as the orbit. This avoids repeated matrix products, while the mixing improvement associated with alternating projections is realized in one step by $G_\infty$.

\subsection{Achieving $G_\infty = \Pi$ via a linear in $n$ number of groups}
Another interesting consequence of the previous section is that in general, taking more groups leads to a decrease in the number of equivalence classes. With sufficient groups, we can obtain a single equivalence class containing the entire state space $\mathcal{X}$. In that case, $G_\infty = \Pi$ and trivially, $G_\infty P G_\infty = \Pi$ as well. 

We now show that it is possible to achieve this with $n-1$ groups, given that $\mathcal{X} = \llbracket n \rrbracket$.\\

\begin{prop} \label{2groupeg}
    For $i \in \llbracket n-1 \rrbracket$, define the two-element group
    $$\mathcal{G}_i = \{e, g_{i+1}\}, \quad g_{i+1} := (1, i+1),$$
    and let $G_i$ be the respective Gibbs kernels. 
    In other words, each $\mathcal{G}_i$ admits a single non-trivial group action that swaps  states 1 and $i+1$. Then the limiting projection $G_\infty = \lim_{m \rightarrow \infty}(G_1\cdots G_{n-1})^m = \Pi$.
\end{prop}

\begin{proof}
    By the construction of $(\mathcal{G}_i)_{i=1}^{n-1}$, for any $x \in \mathcal{X}$, $x \sim 1$. Hence only a single equivalence class exist, and that is equal to $\mathcal{X}$ itself. Then 
    $$G_\infty (x,y) = \frac{\pi(y)}{\sum_{z \in \mathcal{X}} \pi(z)},$$
    and so all its rows are equal to $\pi$. 
\end{proof}

\subsection{Rate of convergence of alternating projections}
The convergence of alternating projections depend greatly on the cosine as defined in \eqref{cosine}. Here we give an upper bound on $c(S_1, S_2)$, that relates closely with the amount of overlapping between the orbit blocks of $S_1$ and $S_2$.\\

\begin{prop} \label{T(j,i)}
    Let $\mathcal{G}_1, \mathcal{G}_2$ be two groups admitting orbits $(\mathcal{O}_i)_{i=1}^{k_1}$ and $(\mathcal{C}_j)_{j=1}^{k_2}$. Define their Gibbs kernel to be $G_1, G_2$, and $S_1, S_2$ to be the respective projection spaces on $\ell^2(\pi)$. Let $T$ be the matrix of size $k_2 \times k_1$, and
    $$T(j,i) := \frac{\pi(\mathcal{O}_i \cap \mathcal{C}_j)}{\sqrt{\pi(\mathcal{O}_i)\pi(\mathcal{C}_j)}}.$$
    
    Then the cosine $c(S_1, S_2) = \sigma_2(T)$, the largest singular value of $T$ less than 1. 

    If all singular values of $T$ are 1, then $c(S_1, S_2)$ is instead $0$, with $G_1G_2 = G_2G_1 = G_\infty$.
\end{prop}

\begin{proof}
    Let $\overline{\pi}_1 = (\pi(\mathcal{O}_1), \dots, \pi(\mathcal{O}_{k_1}))$ and $\overline{\pi}_2 = (\pi(\mathcal{C}_1), \dots, \pi(\mathcal{C}_{k_2}))$ be stationary distributions on the state spaces $(\mathcal{O}_i)_{i=1}^{k_1}$ and $(\mathcal{C}_j)_{j=1}^{k_2}$ respectively. Define the isometries $U: \mathbb{R}^{k_1} \to S_1$ and $V: \mathbb{R}^{k_2} \to S_2$ as per \eqref{U}, and similarly, their adjoints as per \eqref{U*}.

    For any $f \in S = S_1 \cap S_2$, $(G_1G_2 - G_\infty)f = 0$ and so it suffices to restrict our attention to $S^\perp$. On $S^\perp$, $\|G_1 G_2 - G_\infty \|_{op} = \|G_1G_2\|_{op}$. Given that $U$ and $V$ are isometries, and that operator norms are invariant to adjoint (see \cite{Rudin1991}), 
    $$\|V(V^*U)U^*\|_{\ell^2(\pi) \rightarrow \ell^2(\pi)} = \|V^*U\|_{\ell^2(\overline{\pi}_1) \rightarrow \ell^2(\overline{\pi}_2)}.$$

    The linear map $R:= V^*U$ acts on $f \in \ell^2(\overline{\pi}_1)$ by
    $$(Rf)(j) = \sum_{i=1}^{k_1} \frac{\pi(\mathcal{O}_i \cap \mathcal{C}_j)}{\pi(\mathcal{C}_j)}f(i).$$
    
    Writing $\widetilde{f}(i) := \sqrt{\pi(\mathcal{O}_i)} f(i)$,
    \begin{align*}
        \|Rf\|^2_{\ell^2(\pi)} &= \sum_{j=1}^{k_2} \pi(\mathcal{C}_j) \bigg( \sum_{i=1}^{k_1} \frac{\pi(\mathcal{O}_i \cap \mathcal{C}_j)}{\pi(\mathcal{C}_j)} f(i) \bigg)^2 \\
        &= \sum_{j=1}^{k_2} \frac{1}{ \pi(\mathcal{C}_j)} \bigg( \sum_{i=1}^{k_1} \pi(\mathcal{O}_i \cap \mathcal{C}_j) f(i) \bigg)^2\\
        &= \sum_{j=1}^{k_2} (T\widetilde{f})(j)\\
        &= \|Tf\|_2,
    \end{align*}
    where $\|\cdot\|_2$ denotes the norm under the usual Euclidean inner product.

    Since $\|f\|^2_{\ell^2(\overline{\pi}_1)} = \|\widetilde{f}\|_2^2$,
    $$
    \frac{\|Rf\|^2_{\ell^2(\pi)}}{\|f\|^2_{\ell^2(\overline{\pi}_1)}} = \frac{\|T\widetilde{f}\|_2^2}{\|\widetilde{f}\|_2^2}
    $$
    and so $\|R\|_{\ell^2(\overline{\pi}_1) \rightarrow \ell^2(\overline{\pi}_2)} = \|T\|_{2\rightarrow 2}$, the spectral norm of $T$. Equivalently it is also the largest singular value of $T$. 

    Restricting our attention to $S^\perp$, we remove the singular direction associated with $\sigma_1(T)$, and hence $\|G_1G_2 - G_\infty\|_{\ell^2(\pi) \rightarrow \ell^2(\pi)} = \sigma_2(T)$. 

    Now suppose if all singular values of $T$ are 1, and assume that $k_2 \leq k_1$. Then $TT^* = I_{k_2}$. In particular, for $1 \leq j_2 \leq k_2,\ j_1 \neq j_2$,
    $$TT^*(j_1, j_2) = \sum_{i=1}^{k_1} \frac{\pi(\mathcal{O}_i \cap \mathcal{C}_{j_1})\  \pi(\mathcal{O}_i \cap \mathcal{C}_{j_2})}{\pi(\mathcal{O}_i) \sqrt{\pi(\mathcal{C}_{j_1})\pi(\mathcal{C}_{j_2})}} = 0.$$
    This implies that every $\mathcal{O}_i$ must be fully contained within some $\mathcal{C}_j$. Equivalently, $S_2 \subseteq S_1$. Similarly, if $k_1 \leq k_2$ then $S_1 \subseteq S_2$. In any case, $G_1G_2 = G_2G_1 = G_\infty$ and $c(S_1, S_2) = 0$ by \cite{deutsch2001best}, Lemma 9.5.
\end{proof}

As a corollary, we look at the discrete uniform distribution on $\mathcal{X}$ and show that with rqo groups, we can achieve a sizeable convergence rate. \\

\begin{prop} \label{unifeg}
    Let $|\mathcal{X}| = n = mk$, where $m,k$ are both integers, and assume that $\pi$ is the discrete uniform distribution on $\llbracket n \rrbracket$.  Define the groups $\mathcal{G}_1, \mathcal{G}_2$ such that their orbits are given by the partitions
    $$\mathcal{O}_i = \{(i-1)k + 1, \dots,\ ik\},$$
    and 
    $$\mathcal{C}_j = \{j,\ j+m, \dots,\ j+(k-1)m\},$$
    where $i,\ j = 1, \dots, m.$
    
    With this formulation, we can achieve $c(S_1, S_2) \leq m^2/n$ with the Gibbs orbit kernel $G_1, G_2$ satisfying $\lim_{t \rightarrow \infty} (G_1G_2)^t = \Pi$. 

    Suppose $k \geq m$ and $m$ divides $k$, the constructed $G_1, G_2$ can achieve $G_1G_2 = \Pi$.
\end{prop}

\begin{proof}
    Under the uniform distribution,
    $$
    T(j,i) := \frac{\pi(\mathcal{O}_i \cap \mathcal{C}_j)}{\sqrt{\pi(\mathcal{O}_i) \pi(\mathcal{C}_j)}} = \frac{|\mathcal{O}_i \cap \mathcal{C}_j|}{k}.
    $$
    
    Now let $J_m$ be the $m \times m$ matrix of all $1$'s, which is rank 1. Write $T = \frac{1}{m}J_m + A$, and by \cite{Horn_Johnson_1991} Theorem 3.3.16, 
    $$\sigma_2(T) \leq \sigma_2\bigg(\frac{1}{m}J_m\bigg) + \sigma_1(A) = \|A\|_2.$$

    The construction of $\mathcal{O}_i$ and $\mathcal{C}_j$ guarantees that $|\mathcal{O}_i \cap \mathcal{C}_j| = \lfloor k/m \rfloor$ or $\lceil k/m \rceil$, the two integers closest to $k/m$. Hence, for any $i,j$,
    $$|T(j,i) - 1/m| \leq 1/k.$$
    
    Let $\|A\|_1$ and $\|A\|_\infty$ be the maximum absolute column and row sum of $A$. Then we have 
    $$\|A\|_2 \leq \sqrt{\|A\|_1 \cdot \|A\|_\infty} \leq \frac{m}{k} = \frac{m^2}{n},$$
    where one can refer to \cite{Golub2013} Section 2.33 for the first inequality. 

    Furthermore, if $k \geq m$ and $m$ divides $k$, then $|\mathcal{O}_i \cap \mathcal{C}_j|$ is always equals to $k/m$. Then $T = \frac{1}{m}J_m$ and $\sigma_2(T) = 0$. \\
\end{proof}

\begin{corollary}
    With the same context as Proposition \ref{unifeg} above, taking $m = \big\lceil \log n \big\rceil$ gives us 
    $$c(S_1, S_2) \leq \frac{(\log n)^2}{n} = o\bigg( \frac{1}{n^{1-\alpha}}\bigg),$$
    for $\alpha \in (0,1).$
    That is, the convergence rate to $\Pi$ is asymptotically smaller than $n^{\alpha -1}.$
\end{corollary}

 As an application of Proposition \ref{unifeg}, we show that for a state space of size $n = 2^d$, one can achieve $\Pi$ by using $O(d)$ alternating projections. \\
 
\begin{corollary} \label{repeatedalt}
    Let the state space have size $n=2^d$, where $d=2^q$ for some
    positive integer $q$. Then the exact sampler $\Pi$ can be expressed
    using $O(d)=O(\log n)$ alternating projection factors.

    Moreover, the Gibbs orbit kernels at the finest level of the
    construction have orbit blocks of constant size.
\end{corollary}

\begin{proof}
    Let $m = k = \sqrt{n} = 2^{d/2}$, and construct the Gibbs kernels $G_1, G_2$ as described in Proposition \ref{unifeg}. Under this construction, the product $G_1 G_2 = \Pi$ exactly.

    The matrix $G_1$ can be written in block-diagonal form as
    $$
    G_1 = \mathrm{diag}(\Pi_k, \dots, \Pi_k),
    $$
    where each $\Pi_k$ is a $k \times k$ matrix of all entries $1/k$. For each block, one can find a pair of projection matrices $G_1^{(1)}$ and $G_2^{(1)}$ such that 
    $G_1^{(1)} G_2^{(1)} = \Pi_k$, following the same construction.

    Proceeding recursively, at the $r$-th iteration, $G_1^{(r)}$ consists of block-diagonal components of size $2^{d/2^r}$, and each such block can be obtained by applying $2^r$ alternating projection products on smaller sub-blocks. 

    An analogous recursive decomposition applies to $G_2$ (and any subsequent $G_2^{(r)}$), after reordering the indices so that the partitions are expressed in block-diagonal form. Hence, after $r$ recursive levels, the resulting kernels act on disjoint blocks of size $2^{d/2^r}$, and the total number of products required is $2^r$. 

    Taking $r = q = \log_2 d$, we reach the final scale where each block is of constant size, and hence the total number of alternating projection products needed is $O(d) = O(\log n)$.
\end{proof}

To further extend this idea, we propose a model as follows. Suppose we have a state space $\mathcal{X} = \llbracket 0, n-1 \rrbracket$, with $n = 2m^2k$ and $m, k$ as positive integers. Let the stationary distribution on $\mathcal{X}$ be 
$$\pi_\beta(x) \propto e^{\beta |x (\text{mod } 2k) - (k+1)|},$$
which resembles multiple blocks of ``V"-shaped, with a total of $m^2$ modes. 

\begin{figure}[ht]
  \centering
    

\begin{tikzpicture}
\begin{axis}[
  name=pdfaxis,
  width=12cm,
  height=5cm,
  axis lines=middle,
  axis line style={thick,-{Latex[length=3pt]}},
  xlabel={$x$},
  ylabel={$\pi_\beta(x)$},
  xmin=-0.5, xmax=12.5,
  ymin=0, ymax=5.3,
  xtick={0,4,8,12},
  ytick=\empty,
  domain=0:12,
  samples=400,
  thick,
  clip=false,
  enlarge x limits={abs=0.3}
]

\addplot[blue,very thick, smooth]
  ({x},{exp(0.8*abs(mod(x,4)-2))});

\addplot [gray!50, dashed] coordinates {(4,0) (4,5.1)};
\addplot [gray!50, dashed] coordinates {(8,0) (8,5.1)};
\addplot [gray!50, dashed] coordinates {(12,0) (12,5.1)};

\addplot[only marks, mark=*, mark size=1.5pt, blue]
  coordinates {(0,4.9530) (1,2.2255) (2,1.0000) (3,2.2255)};
\addplot[only marks, mark=x, mark size=2pt, blue]
  coordinates {(4,4.9530) (5,2.2255) (6,1.0000) (7,2.2255)};
\addplot[only marks, mark=diamond*, mark size=2pt, blue]
  coordinates {(8,4.9530) (9,2.2255) (10,1.0000) (11,2.2255)};

% --- Store coordinate transformations for later use ---
\pgfplotsextra{
  \path (axis cs:0,0) coordinate (A0);
  \path (axis cs:3,0) coordinate (A3);
  \path (axis cs:4,0) coordinate (A4);
  \path (axis cs:7,0) coordinate (A7);
  \path (axis cs:8,0) coordinate (A8);
  \path (axis cs:11,0) coordinate (A11);
}

\end{axis}

%---------------- Block demarcation arrows ----------------
\begin{scope}
  \draw[latex-latex, thick]
    ([yshift=-17pt]A0) -- ([yshift=-17pt]A3)
    node[midway, below=2pt]{$\mathcal{D}_1$};
  \draw[latex-latex, thick]
    ([yshift=-17pt]A4) -- ([yshift=-17pt]A7)
    node[midway, below=2pt]{$\mathcal{D}_2$};
  \draw[latex-latex, thick]
    ([yshift=-17pt]A8) -- ([yshift=-17pt]A11)
    node[midway, below=2pt]{$\mathcal{D}_3$};
\end{scope}

\end{tikzpicture}
  \caption{Plot of $\pi_\beta(x)$ for $k=2$}
  \label{fig:pdf}
\end{figure}

Partition $\mathcal{X}$ by $(\mathcal{D}_i)_{i=1}^{m^2}$, with $\mathcal{D}_i = \{ x \in  \mathcal{X}: 2(i-1)k \leq x \leq 2ik-1\}$. This formulation ensures that $\pi(\mathcal{D}_i) = 1/m^2$, since each partition has exactly the same points up to cyclic permutation of the indices.  

Now for $i,j \in \llbracket m \rrbracket$, set 
$$\mathcal{O}_i = \bigcup_{l=1}^m \mathcal{D}_{(i-1)m + l} \quad \text{and} \quad \mathcal{C}_j = \bigcup_{l=1}^m \mathcal{D}_{j + (l-1)m}.$$
Here, each orbit $\mathcal{O}_i$ and $\mathcal{C}_j$ consists of exactly $m$ disjoint copies of $\mathcal{D}$'s. For any pair of $i,j$, it is also the case that $\mathcal{O}_i \cap \mathcal{C}_j$ must have exactly 1 such block $\mathcal{D}$. 

By Proposition \ref{T(j,i)}, each entry $T(j,i) = 1/m$ for any $i,j$. Hence, if one constructs the Gibbs orbit kernel $G_1, G_2$ using the partitions $(\mathcal{O}_i)_{i=1}^m$ and $(\mathcal{C}_j)_{j=1}^m$ respectively, $G_1G_2 = \Pi$ must hold again by Proposition \ref{T(j,i)}. 

Finally for $G_1, G_2$, one can then use the technique described in Corollary \ref{repeatedalt}, in which we could then achieve $\Pi$ with $O(\log m)$ projection products. This is a significant improvement over classical dynamics such as Metropolis-Hastings, which in low temperature (e.g. $\beta > 1$ and $m,k$ are chosen such that $k = \Omega(n)$) can have mixing times in the order of $e^n$.

\section{Tuning strategies for choosing $G$} \label{section:tuning}
In Section \ref{section:bestG}, we identified the orbit partition with a fixed number of blocks whose associated Gibbs kernel $G$ minimises
$D_{KL}^\pi(G\|\Pi)$. In practice, however, it is not always feasible to do so, especially if it is not computationally feasible to enumerate through all $\pi(x).$

Recall that our state space is denoted by $\mathcal{X}$ with $|\mathcal{X}| = n$. Let $F: \mathcal{X} \to \mathbb{R}$ be a Hamiltonian function, and
$$\pi_\beta(x) := \frac{1}{Z(\beta)} \exp(-\beta F(x))$$
be the Gibbs distribution associated with the inverse temperature $\beta\geq 0$ with normalisation constant $Z(\beta)$. 

Below we discuss two heuristics for selecting $G$ under $\pi_\beta$, both
aimed at identifying orbit partitions that group states with large stationary
mass.

\subsection{Adaptive heuristic for learning $G$} \label{adaptiveG}

The first heuristic uses the trajectory of a sampler to adaptively learn a
partition that approximates the KL-optimal partition in
Proposition~\ref{bestGbyKL}.

We initialise with $G_0=I$. For a fixed block length $t$ (say 50
iterations), at the $m$-th block we run $t$ iterations using $G_mPG_m$ and
use the states encountered thus far to update the candidate partition.
For a predetermined $k<n$, we identify the $k$ distinct visited states with
the smallest values of the Hamiltonian $F$ and place them in a common
orbit. Since
\[\pi_\beta(x)\propto e^{-\beta F(x)},\]
these states serve as empirical candidates for the high-probability states
that Proposition~\ref{bestGbyKL} places into a single orbit. Every state
outside this orbit is assigned to a singleton orbit. Thus, if
$\mathcal O_1$ denotes the selected orbit of $k$ states, the resulting
partition is
\[\{\mathcal O_1\} \cup \big\{\{x\}:x\in\mathcal X\setminus\mathcal O_1\big\}.\]

Each candidate partition can be realised as the orbit partition of a
finite group action, for example by taking a cycle on the distinguished
orbit and fixing all remaining states. Repeating the procedure periodically
updates the candidate $G$, with the purpose of progressively identifying
the low-energy, high-probability states that should form the distinguished
orbit in Proposition \ref{bestGbyKL}.

\subsection{Initial exploratory chain to learn $G$}
The second heuristic leans towards ``learning" a suitable G in an initial exploration phase using a high-temperature chain. Suppose the goal is to sample from $\pi_{\beta_0}$, with inverse temperature $\beta_0 > 0$ that is potentially large.

At high temperature (small $\beta$), the distribution $\pi_\beta$ tends to be flatter. Usual samplers such as the Metropolis-Hastings algorithm would hence be able to explore the landscape more easily without getting trapped. 

From these empirical frequencies, we construct a partition of the state space by grouping states that appear frequently or are energetically similar, using the same strategy described in Section \ref{adaptiveG} to form the Gibbs kernel $G$. 

This partition determines $G$, which is then fixed and used to form the sampler $GPG$ targeting the actual low-temperature distribution $\pi_{\beta_0}.$

One may also run multiple exploratory chains to learn several Gibbs kernels $G_1, \dots, G_k$, and then apply alternating projections $(G_1 \cdots G_k) P (G_k \cdots G_1)$ on $\pi_{\beta_0}$. The results in Section \ref{section:altproj} apply analogously to this alternating sandwich kernel.

\section*{Acknowledgements}

Michael Choi acknowledges the financial support of the projects A-8001061-00-00, NUSREC-HPC-00001, NUSREC-CLD-00001, A-0000178-01-00, A-0000178-02-00 and A-8003574-00-00 at National University of Singapore.

\newpage
\bibliography{ref}

\begin{thebibliography}{34}
\providecommand{\natexlab}[1]{#1}
\providecommand{\url}[1]{\texttt{#1}}
\expandafter\ifx\csname urlstyle\endcsname\relax
  \providecommand{\doi}[1]{doi: #1}\else
  \providecommand{\doi}{doi: \begingroup \urlstyle{rm}\Url}\fi

\bibitem[Aldous and Fill(2002)]{aldous-fill-2014}
David Aldous and James~Allen Fill.
\newblock {Reversible Markov Chains and Random Walks on Graphs}, 2002.
\newblock Unfinished monograph, recompiled 2014, available at
  \url{https://www.stat.berkeley.edu/users/aldous/RWG/book.html}.

\bibitem[Barker(1965)]{Barker}
AA~Barker.
\newblock {Monte Carlo Calculations of the Radial Distribution Functions for a
  Proton-Electron Plasma}.
\newblock \emph{Australian Journal of Physics}, 18\penalty0 (2):\penalty0
  119--134, 04 1965.
\newblock ISSN 0004-9506.
\newblock \doi{10.1071/PH650119}.
\newblock URL \url{https://doi.org/10.1071/PH650119}.

\bibitem[Birrell et~al.(2022)Birrell, Katsoulakis, Rey-Bellet, and
  Zhu]{birrell2022structurepreservinggans}
Jeremiah Birrell, Markos~A. Katsoulakis, Luc Rey-Bellet, and Wei Zhu.
\newblock {Structure-preserving GANs}, 2022.
\newblock URL \url{https://arxiv.org/abs/2202.01129}.

\bibitem[Bovier and den Hollander(2015)]{Bovier_denHollander_2015}
Anton Bovier and Frank den Hollander.
\newblock \emph{{Metastability: A Potential-Theoretic Approach}}, volume 351.
\newblock Springer International Publishing, Cham, 1st 2015.;1; edition, 2015.
\newblock ISBN 0072-7830.

\bibitem[Br\'emaud(2020)]{Bremaud1999}
Pierre Br\'emaud.
\newblock \emph{{Markov Chains: Gibbs Fields, Monte Carlo Simulation and
  Queues}}, volume 31.;31;.
\newblock Springer International Publishing, Cham, 2nd 2020.;2;2nd 2020;2nd
  2020; edition, 2020.
\newblock ISBN 0939-2475.

\bibitem[Bui et~al.(2012)Bui, Huynh, and
  Riedel]{bui2012automorphismgroupsgraphicalmodels}
Hung~Hai Bui, Tuyen~N. Huynh, and Sebastian Riedel.
\newblock {Automorphism Groups of Graphical Models and Lifted Variational
  Inference}, 2012.
\newblock URL \url{https://arxiv.org/abs/1207.4814}.

\bibitem[Chatterjee and Diaconis(2020)]{Diaconis_2020}
Sourav Chatterjee and Persi Diaconis.
\newblock {Speeding up Markov chains with deterministic jumps}.
\newblock \emph{Probability theory and related fields}, 178\penalty0
  (3-4):\penalty0 1193--1214, 2020.

\bibitem[Chen et~al.(2025)Chen, Katsoulakis, Rey-Bellet, and
  Zhu]{chen2025statisticalguaranteesgroupinvariantgans}
Ziyu Chen, Markos~A. Katsoulakis, Luc Rey-Bellet, and Wei Zhu.
\newblock {Statistical Guarantees of Group-Invariant GANs}, 2025.
\newblock URL \url{https://arxiv.org/abs/2305.13517}.

\bibitem[Choi and Wang(2025)]{choi2025groupaveragedmarkovchainsmixing}
Michael C.~H. Choi and Youjia Wang.
\newblock {Group-averaged Markov chains: mixing improvement}, 2025.
\newblock URL \url{https://arxiv.org/abs/2509.02996}.

\bibitem[Choi and Wolfer(2024)]{Choi_2024}
Michael C.~H. Choi and Geoffrey Wolfer.
\newblock {Systematic Approaches to Generate Reversiblizations of Markov
  Chains}.
\newblock \emph{IEEE Transactions on Information Theory}, 70\penalty0
  (5):\penalty0 3145--3161, 2024.
\newblock \doi{10.1109/TIT.2023.3304685}.

\bibitem[Choi et~al.(2025)Choi, Hird, and Wang]{Choi_Permutation}
Michael C.~H. Choi, Max Hird, and Youjia Wang.
\newblock {Improving the Convergence of Markov Chains via Permutations and
  Projections}.
\newblock \emph{Random Structures \& Algorithms}, 66\penalty0 (4), June 2025.
\newblock ISSN 1098-2418.
\newblock \doi{10.1002/rsa.70016}.
\newblock URL \url{http://dx.doi.org/10.1002/rsa.70016}.

\bibitem[Cohen and
  Welling(2016)]{cohen2016groupequivariantconvolutionalnetworks}
Taco~S. Cohen and Max Welling.
\newblock {Group Equivariant Convolutional Networks}, 2016.
\newblock URL \url{https://arxiv.org/abs/1602.07576}.

\bibitem[Deutsch(2001)]{deutsch2001best}
Frank~R. Deutsch.
\newblock \emph{{Best Approximation in Inner Product Spaces}}.
\newblock Springer New York, New York, NY, 2001.
\newblock ISBN 978-1-4419-2890-0.

\bibitem[Diaconis and Howes(2025)]{Diaconis_2025}
Persi Diaconis and Michael Howes.
\newblock {Random sampling of contingency tables and partitions: Two practical
  examples of the Burnside process}.
\newblock \emph{Statistics and Computing}, 35\penalty0 (6), August 2025.
\newblock ISSN 1573-1375.
\newblock \doi{10.1007/s11222-025-10708-5}.
\newblock URL \url{http://dx.doi.org/10.1007/s11222-025-10708-5}.

\bibitem[Diaconis and Zhong(2021)]{Diaconis_2021_hahn}
Persi Diaconis and Chenyang Zhong.
\newblock {Hahn polynomials and the Burnside process}.
\newblock \emph{The Ramanujan Journal}, 61\penalty0 (2):\penalty0 567–595,
  September 2021.
\newblock ISSN 1572-9303.
\newblock \doi{10.1007/s11139-021-00482-z}.
\newblock URL \url{http://dx.doi.org/10.1007/s11139-021-00482-z}.

\bibitem[Diaconis and Zhong(2025)]{diaconis2025countingnumbergrouporbits}
Persi Diaconis and Chenyang Zhong.
\newblock Counting the number of group orbits by marrying the burnside process
  with importance sampling, 2025.
\newblock URL \url{https://arxiv.org/abs/2501.11731}.

\bibitem[Diaconis et~al.(2025)Diaconis, Lin, and
  Ram]{diaconis2025curiouslyslowlymixingmarkov}
Persi Diaconis, Andrew Lin, and Arun Ram.
\newblock {A curiously slowly mixing Markov chain}, 2025.
\newblock URL \url{https://arxiv.org/abs/2511.01245}.

\bibitem[Frigessi et~al.(1993)Frigessi, Di~Stefano, Hwang, and
  Sheu]{Frigessi1993}
Arnoldo Frigessi, Patrizia Di~Stefano, Chii-Ruey Hwang, and Shuenn-Jyi Sheu.
\newblock {Convergence Rates of the Gibbs Sampler, the Metropolis Algorithm and
  Other Single-Site Updating Dynamics}.
\newblock \emph{Journal of the Royal Statistical Society. Series B,
  Methodological}, 55\penalty0 (1):\penalty0 205--219, 1993.

\bibitem[Ginat(2018)]{ginat2018methodalternatingprojections}
Omer Ginat.
\newblock {The Method of Alternating Projections}, 2018.
\newblock URL \url{https://arxiv.org/abs/1809.05858}.

\bibitem[Golub and Van~Loan(2013)]{Golub2013}
Gene~H. Golub and Charles~F. Van~Loan.
\newblock \emph{Matrix Computations - 4th Edition}.
\newblock Johns Hopkins University Press, Philadelphia, PA, 2013.
\newblock \doi{10.1137/1.9781421407944}.
\newblock URL \url{https://epubs.siam.org/doi/abs/10.1137/1.9781421407944}.

\bibitem[Holley and Stroock(1988)]{Holley_Stroock_1988}
Richard Holley and Daniel Stroock.
\newblock {Simulated annealing via Sobolev inequalities}.
\newblock \emph{Communications in mathematical physics}, 115\penalty0
  (4):\penalty0 553--569, 1988.

\bibitem[Horn and Johnson(1991)]{Horn_Johnson_1991}
Roger~A. Horn and Charles~R. Johnson.
\newblock \emph{{Topics in Matrix Analysis}}.
\newblock Cambridge University Press, 1991.

\bibitem[Jerrum(1993)]{Jerrum_1993}
Mark Jerrum.
\newblock {Uniform sampling modulo a group of symmetries using Markov chain
  simulation}.
\newblock \emph{DIMACS Series in Discrete Mathematics and Theoretical Computer
  Science}, 10:\penalty0 37–47, Jul 1993.
\newblock \doi{10.1090/dimacs/010/04}.

\bibitem[Jerrum et~al.(2004)Jerrum, Son, Tetali, and Vigoda]{Jerrum_2004}
Mark Jerrum, Jung-Bae Son, Prasad Tetali, and Eric Vigoda.
\newblock {Elementary Bounds on Poincaré and Log-Sobolev Constants for
  Decomposable Markov Chains}.
\newblock \emph{The Annals of applied probability}, 14\penalty0 (4):\penalty0
  1741--1765, 2004.

\bibitem[Knuth(1997)]{Knuth_1997}
Donald~E. Knuth.
\newblock \emph{{The Art of Computer Programming, Volume 2: Seminumerical
  Algorithms}}.
\newblock Addison-Wesley, 3rd edition, 1997.
\newblock ISBN 978-0-201-89684-8.
\newblock Fisher--Yates shuffle (Algorithm P).

\bibitem[Kondor and
  Trivedi(2018)]{kondor2018generalizationequivarianceconvolutionneural}
Risi Kondor and Shubhendu Trivedi.
\newblock {On the Generalization of Equivariance and Convolution in Neural
  Networks to the Action of Compact Groups}, 2018.
\newblock URL \url{https://arxiv.org/abs/1802.03690}.

\bibitem[Kou et~al.(2006)Kou, Zhou, and Wong]{Kou_2006}
S.~C. Kou, Qing Zhou, and Wing~Hung Wong.
\newblock {Equi-energy sampler with applications in statistical inference and
  statistical mechanics}.
\newblock \emph{The Annals of Statistics}, 34\penalty0 (4):\penalty0 1581 --
  1619, 2006.
\newblock \doi{10.1214/009053606000000515}.
\newblock URL \url{https://doi.org/10.1214/009053606000000515}.

\bibitem[Levin et~al.(2017)Levin, Peres, Wilmer, Propp, and
  Wilson]{Levin_Peres_Wilmer_2017}
David~A. Levin, Y.~Peres, Elizabeth~L. Wilmer, James Propp, and David~B.
  Wilson.
\newblock \emph{{Markov chains and mixing times}}.
\newblock American Mathematical Society, Providence, Rhode Island, second
  edition, 2017.
\newblock ISBN 1470429624;9781470429621;.

\bibitem[Lim and Choi(2026)]{additive}
Ryan~J. Lim and Michael~C. Choi.
\newblock On additive averaging kernels for finite markov chains.
\newblock \emph{Japanese Journal of Statistics and Data Science}, Aug 2026.
\newblock \doi{10.1007/s42081-026-00373-x}.

\bibitem[Liu(1996)]{Liu1996}
Jun~S. Liu.
\newblock {Metropolized independent sampling with comparisons to rejection
  sampling and importance sampling}.
\newblock \emph{Statistics and computing}, 6\penalty0 (2):\penalty0 113--119,
  1996.

\bibitem[Peskun(1973)]{Peskun1973}
P.~H. Peskun.
\newblock {Optimum Monte-Carlo sampling using Markov chains}.
\newblock \emph{Biometrika}, 60\penalty0 (3):\penalty0 607--612, 12 1973.
\newblock ISSN 0006-3444.
\newblock \doi{10.1093/biomet/60.3.607}.
\newblock URL \url{https://doi.org/10.1093/biomet/60.3.607}.

\bibitem[Rudin(1991)]{Rudin1991}
Walter Rudin.
\newblock \emph{{Functional analysis}}.
\newblock McGraw-Hill, New York, 2nd edition, 1991.
\newblock ISBN 0071009442;0070542368;9780071009447;9780070542365;.

\bibitem[Sherlock(2025)]{sherlock2025reversiblemarkovchainsvariational}
Chris Sherlock.
\newblock {Reversible Markov chains: variational representations and ordering},
  2025.
\newblock URL \url{https://arxiv.org/abs/1809.01903}.

\bibitem[Ying(2022)]{ying2022doubleflipisingmodels}
Lexing Ying.
\newblock {Double Flip Move for Ising Models with Mixed Boundary Conditions},
  2022.
\newblock URL \url{https://arxiv.org/abs/2205.07226}.

\end{thebibliography}

\end{document}